\documentclass[american,british,english]{article}
\usepackage[T1]{fontenc}
\usepackage[latin9]{inputenc}
\usepackage[a4paper]{geometry}
\usepackage{xcolor}
\usepackage{babel}
\usepackage{array}
\usepackage{refstyle}
\usepackage{calc}
\usepackage{mathtools}
\usepackage{dsfont}
\usepackage{amsmath}
\usepackage{lipsum}

\usepackage{amsmath,cite}

\usepackage{amsthm}
\usepackage{amssymb}
\usepackage{graphicx}
\usepackage{wasysym}
\usepackage[all]{xy}
\usepackage[unicode=true,pdfusetitle,
 bookmarks=true,bookmarksnumbered=false,bookmarksopen=false,
 breaklinks=false,pdfborder={0 0 1},backref=false,colorlinks=false]
 {hyperref}

\makeatletter

\AtBeginDocument{\providecommand\defref[1]{\ref{def:#1}}}
\AtBeginDocument{\providecommand\lemref[1]{\ref{lem:#1}}}
\AtBeginDocument{\providecommand\thmref[1]{\ref{thm:#1}}}
\AtBeginDocument{\providecommand\Eqref[1]{\ref{Eq:#1}}}
\AtBeginDocument{\providecommand\assuref[1]{\ref{assu:#1}}}
\AtBeginDocument{\providecommand\remref[1]{\ref{rem:#1}}}
\AtBeginDocument{\providecommand\corref[1]{\ref{cor:#1}}}
\AtBeginDocument{\providecommand\Figref[1]{\ref{Fig:#1}}}
\providecommand{\tabularnewline}{\\}
\RS@ifundefined{subsecref}
  {\newref{subsec}{name = \RSsectxt}}
  {}
\RS@ifundefined{thmref}
  {\def\RSthmtxt{theorem~}\newref{thm}{name = \RSthmtxt}}
  {}
\RS@ifundefined{lemref}
  {\def\RSlemtxt{lemma~}\newref{lem}{name = \RSlemtxt}}
  {}

\theoremstyle{plain}
\newtheorem{thm}{\protect\theoremname}
\theoremstyle{definition}
\newtheorem{defn}[thm]{\protect\definitionname}
\theoremstyle{remark}
\newtheorem{rem}[thm]{\protect\remarkname}
\theoremstyle{plain}
\newtheorem{assumption}[thm]{\protect\assumptionname}
\theoremstyle{plain}
\newtheorem{lem}[thm]{\protect\lemmaname}
\theoremstyle{plain}
\newtheorem{cor}[thm]{\protect\corollaryname}
\theoremstyle{remark}
\newtheorem*{notation*}{\protect\notationname}
\theoremstyle{remark}
\newtheorem*{acknowledgement*}{\protect\acknowledgementname}

\@ifundefined{date}{}{\date{}}
\usepackage{dsfont}
\usepackage[nottoc,numbib]{tocbibind}
\usepackage{xypic}
\usepackage{bbm}

\newref{rem}{name=Remark~}
\newref{def}{name=Definition~}
\newref{lem}{name=Lemma~}
\newref{thm}{name=Theorem~}
\newref{prop}{name=Proposition~}
\newref{cor}{name=Corollary~}
\newref{assu}{name=Assumption~}
\newref{cri}{name=Criterion~}

\makeatletter
\newcommand{\xyR}[1]{%
\makeatletter
\xydef@\xymatrixrowsep@{#1}
\makeatother} % end of \xyR

\makeatletter
\newcommand{\xyC}[1]{%
\makeatletter
\xydef@\xymatrixcolsep@{#1}
\makeatother} % end of \xyC

\DeclareMathOperator*{\esssup}{ess\,sup}

\usepackage{fancyhdr}
\makeatother

\newcommand\blfootnote[1]{%
  \begingroup
  \renewcommand\thefootnote{}\footnote{#1}%
  \addtocounter{footnote}{-1}%
  \endgroup
}

\addto\captionsamerican{\renewcommand{\acknowledgementname}{Acknowledgement}}
\addto\captionsamerican{\renewcommand{\assumptionname}{Assumption}}
\addto\captionsamerican{\renewcommand{\corollaryname}{Corollary}}
\addto\captionsamerican{\renewcommand{\definitionname}{Definition}}
\addto\captionsamerican{\renewcommand{\lemmaname}{Lemma}}
\addto\captionsamerican{\renewcommand{\notationname}{Notation}}
\addto\captionsamerican{\renewcommand{\remarkname}{Remark}}
\addto\captionsamerican{\renewcommand{\theoremname}{Theorem}}
\addto\captionsbritish{\renewcommand{\acknowledgementname}{Acknowledgement}}
\addto\captionsbritish{\renewcommand{\assumptionname}{Assumption}}
\addto\captionsbritish{\renewcommand{\corollaryname}{Corollary}}
\addto\captionsbritish{\renewcommand{\definitionname}{Definition}}
\addto\captionsbritish{\renewcommand{\lemmaname}{Lemma}}
\addto\captionsbritish{\renewcommand{\notationname}{Notation}}
\addto\captionsbritish{\renewcommand{\remarkname}{Remark}}
\addto\captionsbritish{\renewcommand{\theoremname}{Theorem}}
\addto\captionsenglish{\renewcommand{\acknowledgementname}{Acknowledgement}}
\addto\captionsenglish{\renewcommand{\assumptionname}{Assumption}}
\addto\captionsenglish{\renewcommand{\corollaryname}{Corollary}}
\addto\captionsenglish{\renewcommand{\definitionname}{Definition}}
\addto\captionsenglish{\renewcommand{\lemmaname}{Lemma}}
\addto\captionsenglish{\renewcommand{\notationname}{Notation}}
\addto\captionsenglish{\renewcommand{\remarkname}{Remark}}
\addto\captionsenglish{\renewcommand{\theoremname}{Theorem}}
\providecommand{\acknowledgementname}{Acknowledgement}
\providecommand{\assumptionname}{Assumption}
\providecommand{\corollaryname}{Corollary}
\providecommand{\definitionname}{Definition}
\providecommand{\lemmaname}{Lemma}
\providecommand{\notationname}{Notation}
\providecommand{\remarkname}{Remark}
\providecommand{\theoremname}{Theorem}

\begin{document}
\global\long\def\X{\mathbb{X}}%
\global\long\def\x{\mathbbm x}%
\global\long\def\xx{\mathbf{x}}%
\global\long\def\xn{\x^{(n)}}%
\global\long\def\th{\mathtt{F}_{n}}%
\global\long\def\Fth{\mathtt{F}}%
\global\long\def\O{\mathcal{Q}}%
\global\long\def\y{\mathbbm y}%

\global\long\def\B{\mathcal{B}}%
\global\long\def\E{\mathbb{E}}%
\global\long\def\P{\mathbb{P}}%

\global\long\def\cC{\mathcal{C}}%
\global\long\def\cT{\mathcal{T}}%
\global\long\def\cU{\mathcal{U}}%
\global\long\def\cV{\mathcal{V}}%
\global\long\def\cN{\mathcal{N}}%

\global\long\def\S{\mathcal{S}}%
\global\long\def\R{\mathbb{R}}%
\global\long\def\Rd{\mathbb{\R}^{d}}%
\global\long\def\N{\mathbb{N}}%
\global\long\def\Q{\mathbb{Q}}%
\global\long\def\Z{\mathbb{Z}}%
\global\long\def\F{\mathcal{F}}%
\global\long\def\H{\mathcal{H}}%
\global\long\def\del{\partial}%
\global\long\def\M{\mathcal{M}}%
\global\long\def\A{\mathcal{A}}%
\global\long\def\p{\mathcal{P}}%
\global\long\def\mup{\mu_{\p}}%

\global\long\def\dx{\text{\,\ensuremath{\mathrm{d}}}}%
\global\long\def\scale{\text{\ensuremath{\mathrm{scale}}}}%
\global\long\def\poi{\mathrm{poi}}%

\noindent 
\global\long\def\I{\mathds{1}}%
\global\long\def\L{\mathcal{L}}%
\global\long\def\K{\mathcal{K}}%
\global\long\def\k{\kappa}%
\global\long\def\eps{\varepsilon}%
\global\long\def\gc{\mathbf{G}^{\complement}}%
\global\long\def\bG{\mathbf{G}}%
\global\long\def\tauxe{\tau_{\frac{x}{\eps}}\omega}%
\global\long\def\bign{\big\|}%
\global\long\def\biga{\big|}%
\global\long\def\bm{\boldsymbol{\boxminus}}%
\global\long\def\bigl{\big\langle}%
\global\long\def\bigr{\big\rangle}%

\global\long\def\d{\mathrm{d}}%
\global\long\def\rmD{\mathrm{D}}%
\global\long\def\dist{\mathrm{dist}}%
\global\long\def\loc{\mathrm{loc}}%
\global\long\def\pot{\mathrm{pot}}%
\global\long\def\sol{\mathrm{sol}}%
\global\long\def\ueo{u_{\omega}^{\eps}}%
\global\long\def\uexn{u^{\eps}}%
\global\long\def\Qexn{Q_{\xn}^{\eps}}%

\noindent 
\global\long\def\Ball#1#2{\mathbb{B}_{#1}{\left(#2\right)}}%
\global\long\def\of#1{{\left(#1\right)}}%
\global\long\def\norm#1{\left\Vert #1\right\Vert }%

\global\long\def\weakto{\rightharpoonup}%

\global\long\def\mugammapalm{\mu_{\Gamma,\mathcal{P}}}%
\global\long\def\mupalm{\mu_{\mathcal{P}}}%

\global\long\def\muomega{\mu_{\omega}}%
\global\long\def\muomegaeps{\mu_{\omega}^{\eps}}%
\global\long\def\mugammaomega{\mu_{\Gamma(\omega)}}%
\global\long\def\mugammaomegaeps{\mu_{\Gamma(\omega)}^{\eps}}%

\global\long\def\ue{u^{\eps}}%

\global\long\def\tsto{\stackrel{2s}{\to}}%
\global\long\def\tsweakto{\stackrel{2s}{\weakto}}%

%\author{Alexander Hinsen}
%\address[Alexander Hinsen]{Weierstrass Institute Berlin, Mohrenstr. 39, 10117 Berlin, Germany}
%\email{Alexander.Hinsen@wias-berlin.de}
%
%\keywords{Interacting particle systems; random graphs; survival; extinction; percolation; Boolean model}
%\subjclass[2010]{Primary 60J25; secondary 60K35, 60K37}

%\author{
% Martin Heida\footnotemark[1]\footnotetext[1]{TEST}\\
%  \texttt{Martin.Heida@wias-berlin.de}
%\and
% Martin Heida\footnotemark[1]\\
%  \texttt{Martin.Heida@wias-berlin.de}
%  \and
%   Martin Heida\footnotemark[1]\\
%  \texttt{Martin.Heida@wias-berlin.de}
%}
%\thanks{TEST}

\author{Martin Heida\thanks{Weierstrass Institute Berlin, Mohrenstr. 39, 10117 Berlin, Germany,  \texttt{Martin.Heida@wias-berlin.de}}
\and
Benedikt Jahnel\thanks{Weierstrass Institute Berlin, Mohrenstr. 39, 10117 Berlin, Germany,  \texttt{Benedikt.Jahnel@wias-berlin.de}}
\and
Anh Duc Vu\thanks{Weierstrass Institute Berlin, Mohrenstr. 39, 10117 Berlin, Germany,  \texttt{AnhDuc.Vu@wias-berlin.de}}}

\title{Stochastic Homogenization on Irregularly Perforated Domains}
\maketitle
\begin{abstract}
We study stochastic homogenization of a quasilinear parabolic PDE
with nonlinear microscopic Robin conditions on a perforated domain.
The focus of our work lies on the underlying geometry that does not
allow standard homogenization techniques to be applied directly. Instead
we prove homogenization on a regularized geometry and demonstrate
afterwards that the form of the homogenized equation is independent
from the regularization. Then we pass to the regularization limit
to obtain the anticipated limit equation. Furthermore, we show that
Boolean models of Poisson point processes are covered by our approach.
\end{abstract}
\tableofcontents{}

\section*{Introduction}

Soon after the groundbreaking introduction of stochastic homogenization
by \foreignlanguage{american}{Papanicolaou} and Varadhan \cite{papanicolaou1979boundary}
and Kozlov \cite{kozlov1979averaging}, research developed a natural
interest in the homogenization on randomly perforated domains. A good
summary over the existing methods up to 1994 can be found in \cite{kozlov1994homogenization}.
By the same time, Zhikov \cite{zhikov1993averaging} provided a homogenization
result for linear parabolic equations on stationary randomly perforated
domains. It then became silent for a decade. In \cite{Zhikov2006},
Zhikov and Piatnitsky reopened the case by introducing the stochastic
two-scale convergence as a generalization of \cite{Nguetseng1989,allaire1992homogenization,Zhikov2000}
to the stochastic setting, particularly to random measures that comprise
random perforations and random lower-dimensional structures in a natural
way. The method was generalized to various applications in discrete
and continuous homogenization \cite{MP07,Faggionato2008,flegel2019homogenization}
and recently also to an unfolding method \cite{neukammvarga2018stochastic,heidaneukamm2021stochastic}.

Concerning the homogenization on randomly perforated domains, there
seem to be few results in the literature, with \cite{guillen2015quasistatic,franchi2019mathematical,piatnitski2020homogenization}
being the closest related work from the PDE point of view. We emphasize
that there is a further discipline in stochastic homogenization, studying
critical regimes of scaling for holes in a perforated domain of the
stokes equation, see \cite{giunti2020convergence} and references
therein. 

In this work, we focus on the geometric aspects in the homogenization
of quasilinear parabolic equations and go beyond any recent assumptions
on the random geometry. Given $\eps>0$, we consider for a bounded
domain $Q\subset\Rd$ perforated by a random set $G^{\eps}$ and write
$Q^{\eps}:=Q\backslash G^{\eps}$. Typically, $G^{\eps}\approx\eps G$
where $G$ is a stationary random set and $G^{\eps}$ is additionally
regularized close to $\partial Q$ \cite{guillen2015quasistatic,franchi2019mathematical,piatnitski2020homogenization}.
We then study the following PDE on $Q^{\eps}$ for the time interval
$I=[0,T]$:

\begin{equation}
\begin{aligned}\del_{t}u^{\eps}-\nabla\cdot\left(A(u^{\eps})\,\nabla u^{\eps}\right) & =f &  & \text{in }I\times Q^{\eps}\\
A(u^{\eps})\,\nabla u^{\eps}\cdot\nu & =0 &  & \text{on }I\times\del Q\\
A(u^{\eps})\,\nabla u^{\eps}\cdot\nu & =h(u^{\eps}) &  & \text{on }I\times\del Q^{\eps}\backslash\del Q\\
u^{\eps}(0,x) & =u_{0}(x) &  & \text{in }Q^{\eps}\,.
\end{aligned}
\label{eq:System-in-Eps}
\end{equation}
In case of a fully linear PDE, i.e. $h(\cdot)=const$ and $A(\cdot)=const$,
this problem was homogenized already in the aforementioned \cite{zhikov1993averaging}.
In this linear case one benefits from the regularity of the limit
solution and the weak convergence of the $\eps$-solutions that is
given a priori.

The nonlinear case is, however, more difficult. Weak convergence of
solutions is no longer enough and one needs to establish strong convergence
of the $u^{\eps}$. Typical assumptions in the literature, such as
minimal smoothness (see \defref{minimal-smooth}) of $G$ and uniform
boundedness of the holes, ensure the existence of uniformly bounded
extension operators $\text{\ensuremath{\cU}}_{\eps,\bullet}:W^{1,2}(Q^{\eps})\to W^{1,2}(Q)$
(see, \cite{guillen2015quasistatic}). This in turn implies weak compactness
of $\cU_{\eps,\bullet}u^{\eps}$ in $W^{1,2}(Q)$, a property of uttermost
importance to pass to the homogenization limit in the nonlinear terms.
Other approaches are thinkable, e.g. exploiting the Frechet--Riesz--Kolmogorov
compactness theorem, but in application the prerequisites are hard
to prove.

If all limit passages go through, the homogenized limit as $\eps\to0$
reads for some positive definite matrix $\A_{(G)}$ as 
\begin{align}
C_{1,(G)}\del_{t}u-\text{div}\left(A(u)\A_{(G)}\,\nabla u\right)-C_{2,(G)}h(u) & =C_{1,(G)}f &  & \text{in }I\times Q\nonumber \\
A(u)\A_{(G)}\,\nabla u\cdot\nu & =0 &  & \text{on }I\times\del Q\label{eq:System-Homogenized-G}\\
u(0,x) & =C_{1,(G)}u_{0}(x) &  & \text{in }Q\,,\nonumber 
\end{align}
which represents the macroscopic \foreignlanguage{american}{behavior}
of our object. We note at this point that positivity of $\A_{(G)}$
is in general non-trivial but can be shown for minimally smooth domains
\cite{guillen2015quasistatic} and other examples (see Sections \ref{sec:Sufficient-Condition-for-Effective-Conductivity}
and \ref{sec:PPP}).

\noindent \medskip{}

\noindent Unfortunately, canonical perforation models are neither
minimally smooth nor do they come up with uniformly bounded holes.
Our toy model of choice will be the Boolean model $\Xi\X_{\poi}:=\overline{\cup_{x\in\X_{\poi}}\Ball rx}$
(see \defref{boolean-model-filled}) driven by a Poisson point process
$\X_{\poi}$. It clearly reveals the following general issues for
the homogenization analysis: 
\begin{enumerate}
\item $\Xi\X_{\poi}{}^{\complement}=\R^{d}\backslash\Xi\X_{\poi}$ is not
connected. This happens due to areas that are encircled.
\item Two distinct balls can lie arbitrarily close to each other or --
in case they intersect -- have arbitrary small overlap. This implies
that
\begin{itemize}
\item the connected components in $\Xi\X_{\poi}$ develop arbitrarily large
local Lipschitz constants: Two balls of equal radius intersecting
at an angle $\alpha$ have the Lipschitz constant $\tan((\pi-\alpha)/2)$
at the points of intersection, and
\item there is no $\delta>0$ such that for every $p\in\partial\Xi\X_{\poi}{}^{\complement}$
the surface $\Ball{\delta}p\cap\partial\Xi\X_{\poi}{}^{\complement}$
is a graph of a function: If $x,y\in\X_{\poi}$ with $|x-y|=2r+\eta$
and $\left|p-x\right|=r$, $\left|p-y\right|=r+\eta$, $\Ball{\delta}p\cap\partial\Xi\X_{\poi}$
can be a graph only if $\delta<\eta$.
\end{itemize}
\end{enumerate}
The first issue can be fixed by considering a ``filled-up model''
$\bm\X_{\poi}$ in \defref{boolean-model-filled}. The second issue
poses an actual problem though. In a recent work \cite{heida2021stochastic},
one of the authors has shown that in some cases an extension operator
$\text{\ensuremath{\cU}}_{\eps,\bullet}:W^{1,p}(Q^{\eps})\to W^{1,r}(Q)$,
$1\leq r<p$, can be constructed for some geometries including the
Boolean model (strictly speaking this was shown for an extension from
the balls to the complement in the percolation case). However, \cite{heida2021stochastic}
also suggests that the Boolean model for the Poisson point process
requires $p>2$ for $\text{\ensuremath{\cU}}_{\eps,\bullet}$ to be
properly defined for some $r>0$. 

Due to these severe analytical difficulties, we are in need to try
other approaches to the problem. Our approach includes the following
steps:
\begin{enumerate}
\item Given a general stationary ergodic (admissible) random point process
$\X$, we construct a regularization $\X^{(n)}:=\th\X$ (see \defref{Thinning-Maps})
such that the set $\bm\X^{(n)}$ is uniformly minimally smooth for
given $n\in\N$. 
\item Given $n\in\N$, we perform homogenization for the smoothed geometry
$\bm\X^{(n)}$ instead of $\bm\X$ (see \lemref{Existence-of-Solution-eps}). 
\item We pass to the limit $n\to\infty$ to obtain the anticipated homogenized
limit problem (see \thmref{Main-Theorem}). This happens under the
assumption that $\bm\X{}^{\complement}$ is statistically connected
(see \defref{Statistically-Connected}).
\item We show that the Poisson point process in the subcritical regime is
a valid example for our general homogenization result (see Section
\ref{sec:PPP}).
\end{enumerate}
We are thus in a position to prove an indirect homogenization result.
This seems to us an appropriate intermediate step on the way to a
full homogenization result, which may be achieved in the future using
further developed homogenization techniques based on a better understanding
of the interaction of geometry and homogenization.

\medskip{}
This paper is structured in the following way:
\begin{itemize}
\item In Section \ref{sec:Setting-Definitions-Main-Result}, we introduce
the core objects and state the main result. This includes the thinned
point processes $\X^{(n)}$ and its filled-up Boolean model $\bm\X^{(n)}$.
\item In Section \ref{sec:Properties-of-Thinning-Map}, we prove relevant
properties of the thinning map and the thinned point processes, most
importantly minimal smoothness of $\bm\X^{(n)}$ (\thmref{Minimal-Smoothness-By-F-n})
and $\bm\X^{(n)}\to\bm\X$ in a certain sense (\thmref{Approximation-Properties}).
\item Section \ref{sec:Convergences-of-Cell-Solutions} deals with the cell
solutions and the definition of the effective conductivity $\A$.
\item The homogenization theory for minimally smooth holes is sketched in
Section \ref{sec:Proof-of-Lemma} on the basis of stochastic two-scale
convergence. Due to the considerations in Section \ref{sec:Convergences-of-Cell-Solutions},
the underlying probability space is a compact separable metric space.
\item In Section \ref{sec:Main-Theorem}, we show that the homogenized solutions
to \Eqref{System-Homogenized-G} for $G=\bm\X^{(n)}$ converge and
that their limit is a solution to the anticipated limit problem for
$G=\bm\X$.
\item Section \ref{sec:Sufficient-Condition-for-Effective-Conductivity}
establishes a criterion for statistical connectedness (non-degeneracy
of the effective conductivity $\A$) using percolation channels. This
follows the ideas in \cite[Chapter 9]{kozlov1994homogenization} where
a discrete model was considered.
\item In Section \ref{sec:PPP}, we show that the Poisson point process
$\X_{\poi}$ is indeed admissible which follows from readily available
percolation results. Showing statistical connectedness of $\bm\X_{\poi}^{\complement}$
is much harder. We do so using the criterion established in Section
\ref{sec:Sufficient-Condition-for-Effective-Conductivity} and a version
of \cite[Theorem 11.1]{kesten1982percolation}. As the original \cite[Theorem 11.1]{kesten1982percolation}
is a statement about percolation channels on the $\Z^{2}$-lattice,
we need to fit the statement and proof to our setting.
\end{itemize}

\section*{Notation}

\noindent\fbox{\begin{minipage}[t]{1\columnwidth - 2\fboxsep - 2\fboxrule}%
\textbf{General notation}
\begin{itemize}
\item $\mathcal{M}(\R^{d})$: Space of Radon measures on $\R^{d}$ equipped
with the vague topology
\item $\S(\R^{d})\subset\M(\R^{d})$: Space of boundedly finite point clouds/point
measures in $\R^{d}$
\item $A^{\complement}$: Complement of a set $A$
\item $\B(X)$: Borel-$\sigma$-algebra of the topological space $X$
\item $\L^{d}$: $d$-dimensional Lebesgue-measure
\item $\H^{d}$: $d$-dimensional Hausdorff-measure
\item $\H_{\llcorner A}^{d}$: Restriction of $\H^{d}$ to $A$, i.e. $\H_{\llcorner A}^{d}(B):=\H^{d}(B\cap A)$
\item $o:=0_{\R^{d}}\in\R^{d}$: Origin in $\R^{d}$
\item $\I_{A}$: Indicator/characteristic function of a set A
\end{itemize}
\end{minipage}}

\medskip{}
\noindent\fbox{\begin{minipage}[t]{1\columnwidth - 2\fboxsep - 2\fboxrule}%
\textbf{Specific notation introduced later}
\begin{itemize}
\item $\Ball rA$: Open $r$-neighborhood around $A$. (\defref{boolean-model-filled})
\item $\Xi\x$ and $\bm\x$: Boolean model of $\x$ and its filled version
(\defref{boolean-model-filled})
\item $\cC_{\x}(x)$: Cluster of $x$ in $\x\in\S(\R^{d})$ (\defref{Thinning-Maps})
\item $\x^{(n)}$ for $\x\in\S(\R^{d})$: $\x^{(n)}=\th\x$ with thinning
map $\th$ (\defref{Thinning-Maps})
\item $Q_{\x}^{\eps}$ and $J^{\eps}(Q,\x)$: Perforated domain and index
set generating perforations (\defref{perforation-Q})
\item $\tau_{x}:\ \M(\R^{d})\to\M(\R^{d})$: Shift-operator in $\M(\R^{d})$
(\defref{Random-Measure})
\item $\lambda(\mu)$: Intensity of random measure $\mu$ (\defref{Random-Measure})
\item $\mu_{\x}$: $\H^{d-1}$ restricted to $\del\bm\x$ (\defref{Surface-Measure-of-BM})
\item $\A$ and $\alpha_{\A}$: Effective conductivity and smallest eigenvalue
of $\A$ (\defref{Effective-Conductivity})
\item $\cU$ and $\cT$: Extension and trace operators (\thmref{extension-trace-operators}
and \thmref{Extension-Trace-Estimates-eps})
\item $\mu^{\eps}$: Scaled measure (\assuref{TS-assumption})
\end{itemize}
\end{minipage}}

\section{Setting and main result \label{sec:Setting-Definitions-Main-Result}}

\subsection{Generating minimally smooth perforations}

We start by introducing some concepts from the theory of point processes.
We will not formulate the concepts in full generality but only as
general as needed for our purpose. Let $d\geq2$ and let $\S(\R^{d})$
be the space of boundedly finite point clouds in $\R^{d}$ (i.e. point
clouds without accumulation points) with the Fell topology and $\M(\Rd)$
the space of Radon measures with the vague topology. Every $\x\in\S(\Rd)$
can be identified with a Borel measure through the measurable correspondence
\[
\x(A)=\sum_{x\in\x}\delta_{x}(A)\,.
\]
Hence we identify $\S(\R^{d})\subset\M(\R^{d})$. 

\begin{figure}
\centering{}\includegraphics[width=4cm]{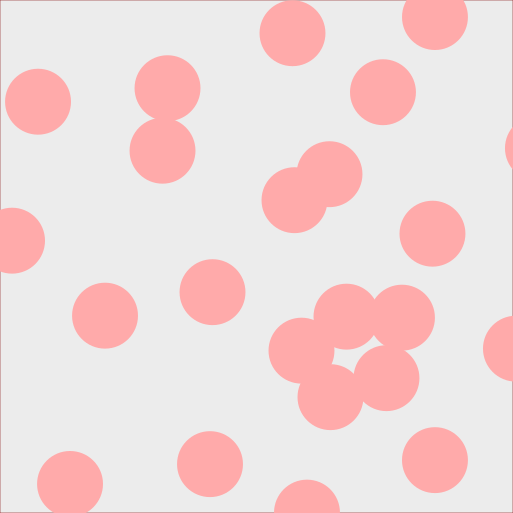}\hspace{0.1\textwidth}\includegraphics[width=4cm]{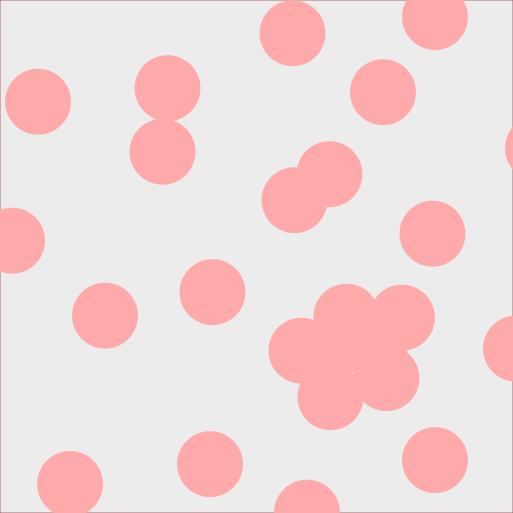}\caption{\label{fig:Filling}Initial Boolean model $\Xi\protect\x$ vs filled-up
Boolean model $\protect\bm\protect\x$.}
\end{figure}

\medskip{}
Our perforation model of interest is the Boolean model driven by a
point cloud $\x\in\S(\R^{d})$. While it is a natural way to generate
perforations, we need to fill it up so that its complement is connected
for suitable $\x$.
\begin{defn}[Boolean model $\Xi$ of a point cloud and filled-up model $\bm$ (see
Figure \ref{fig:Filling})]
\label{def:boolean-model-filled}\ \\
Let $\x\in\S(\R^{d})$. The \emph{Boolean model} of $\x$ for a radius
$r>0$ is
\[
\Xi\x:=\overline{\bigcup_{x\in\x}\Ball rx}=\overline{\Ball r{\x}}\,,
\]
where $\Ball rx$ is the open ball of radius $r$ around $x$ and
$\Ball rA:=\cup_{x\in A}\Ball rx$. \\
We define the \emph{filled-up Boolean model} $\bm\x$ of $\x$ for
radius $r$ through its complement, i.e.
\begin{align*}
\bm\x{}^{\complement}:=\big\{ x\in\R^{d}\,\vert\, & \exists\gamma:\,[0,\infty)\to\Xi\x{}^{\complement}\text{ continuous and }\gamma(0)=x,\,\limsup_{t\to\infty}\biga\gamma(t)\biga=\infty\big\}.
\end{align*}
\end{defn}

\begin{rem}
\label{rem:drawback-boolean}We observe that 
\[
\Xi(\x+x)=\Xi(\x)+x\,,\qquad\bm(\x+x)=\bm(\x)+x\,.
\]
As discussed in the introduction, we need to ``smoothen'' the geometry
in order to be able to apply standard homogenization methods. Given
a Lipschitz domain $P\subset\Rd$, we define for $p\in\partial P$
\[
\delta{\left(p\right)}:=\frac{1}{2}\sup_{\delta'>0}\left\{ \,\partial P\text{ is }\text{Lipschitz-graph in }\Ball{\delta'}p\right\} \,,
\]
and because $\delta:\,\partial P\to\R_{\geq0}$ is continuous \cite{heida2021stochastic},
we can define for bounded $P$
\[
\delta(P):=\min_{p\in\partial P}\delta(p)\,.
\]
\end{rem}

\begin{figure}
\centering{}\includegraphics[width=4cm]{initial_process}$\qquad$\includegraphics[width=4cm]{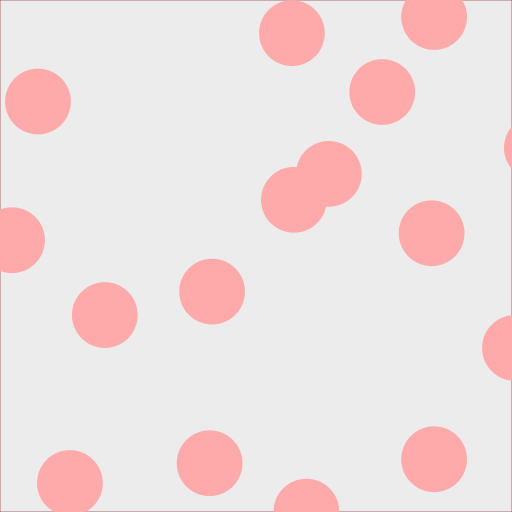}$\qquad$\includegraphics[width=4cm]{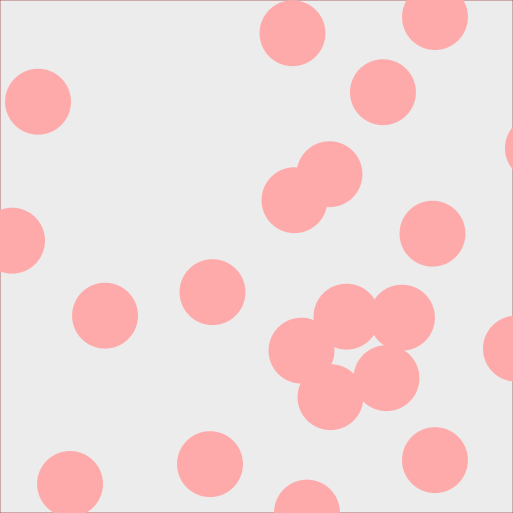}\caption{\label{fig:Thinning}Thinning of point clouds under $\protect\th$
pictured via the Boolean model $\Xi$. From left to right are $\protect\x$,
$\protect\x^{(2)}$ and $\protect\x^{(5)}$.}
\end{figure}

\begin{defn}[Thinning maps $\th$(see Figure \ref{fig:Thinning})]
\label{def:Thinning-Maps}\ \\
Let $\x\in\S(\R^{d})$ be a point cloud. We denote the \emph{cluster}
of $x$ in $\x$ by 
\[
\cC_{\x}(x):=\{y\in\x\,\vert\,\exists\text{path from }x\text{ to }y\text{ inside }\Xi\x\}.
\]
We set
\begin{align*}
\Fth_{1,n}\x & :=\left\{ x\in\x\,\vert\,\forall y\in\x:\,d(x,y)\notin\big(0,\frac{1}{n}\big)\cup\big(2r-\frac{1}{n},\,2r+\frac{1}{n}\big)\right\} \\
\Fth_{2,n}\x & :=\left\{ x\in\x\,\vert\,\#\mathcal{C}_{\x}(x)\leq n,\,\delta\of{\Ball r{\cC_{\x}(x)}}\geq\frac{1}{n}\right\} 
\end{align*}
and define the \emph{thinning map} $\th$ 
\[
\th:\ \S(\R^{d})\to\S(\R^{d})\,,\qquad\xn:=\th\x:=\big(\Fth_{2,n}\circ\Fth_{1,n}\big)\x\,.
\]
\end{defn}

$\th$ can be understood as a generalization of the classical Matérn
construction \cite{matern1986spatial,stoyan1987stochgeo}. For an
arbitrary $\x\in\S(\R^{d})$, we have that $(\bm\xn)^{\complement}$
is always minimally smooth (see \defref{minimal-smooth}). Furthermore,
if $\X$ is a stationary point process (as defined later), then the
same holds for $\X^{(n)}=\th\X$. We note that $\th$ is in general
not monotone in $n$, i.e. $\Fth_{m}\x\not\subset\th\x$ for $m\leq n$.

\medskip{}
Given a scale $\eps>0$, we define the perforation domain $Q^{\eps}$
such that the perforations have some minimal distance from the boundary
$\del Q$:
\begin{defn}[Perforation of domain $Q^{\eps}$]
\label{def:perforation-Q}\ \\
Let $\x\in\S(\Rd)$. We set
\begin{align*}
J^{\eps}(\x,Q) & :=\left\{ x\in\x\,\vert\,\dist\of{\eps\,\cC_{\x}(x)\,,\,Q^{\complement}}>2\eps r\right\} ,\qquad G_{\x}^{\eps}:=\eps\bm(J^{\eps}(\x,Q))
\end{align*}
as well as the \emph{perforated domain}
\[
Q_{\x}^{\eps}:=Q\,\backslash\,G_{\x}^{\eps}\,.
\]
One quickly verifies that $Q_{\xn}^{\eps}$ is minimally smooth (\defref{minimal-smooth}),
see \thmref{Minimal-Smoothness-By-F-n}.
\end{defn}

\subsection{Homogenization for minimally smooth perforations}

\noindent We make the following parameter assumptions on our partial
differential equation (\Eqref{System-in-Eps}).
\begin{assumption}[Parameters of PDE]
\label{assu:PDE-Parameters}\ \\
Let $I=[0,T]\subset\R$ and $Q\subset\R^{d}$ be a bounded, connected
open domain. We assume that
\begin{itemize}
\item $u_{0}\in W^{1,2}(Q)$
\item $f\in L^{2}(I;L^{2}(Q))$
\item $h:\ \R\to\R$ is Lipschitz continuous with Lipschitz constant $L_{h}$
\item $A:\ \R\to\R$ is continuous with $0<\inf(A)$ and $\sup(A)<\infty$.
\end{itemize}
Generalized time derivatives will always be considered under the evolution
triple $W^{1,2}(Q)\hookrightarrow L^{2}(Q)\hookrightarrow W^{1,2}(Q)^{*}$
or $W^{1,2}(\Qexn)\hookrightarrow L^{2}(\Qexn)\hookrightarrow W^{1,2}(\Qexn)^{*}$
in the case of a perforated domain $\Qexn$.
\end{assumption}

\begin{lem}[Solution to PDE for minimally smooth holes]
\label{lem:Existence-of-Solution-eps}\ \\
Let $\x\in\S(\R^{d})$ and $n\in\N$. Under \assuref{PDE-Parameters},we
have on $\Qexn$:\\
There exists a weak solution $\uexn\in L^{2}(I;W^{1,2}(\Qexn))$ with
generalized time derivative $\del_{t}\uexn\in L^{2}(I;W^{1,2}(\Qexn)^{*})$
to \Eqref{System-in-Eps}.\\
This $u^{\eps}$ satisfies for some $C>0$ depending only on $Q$,
$n$, $f$ and $u_{0}$ but not on $\eps$ 
\begin{align*}
\esssup_{t\in I}\bign\uexn(t)\bign_{L^{2}(\Qexn)}+\bign\uexn\bign_{L^{2}(I;W^{1,2}(\Qexn))} & \leq C\\
\bign\del_{t}\uexn\bign_{L^{2}(I;W^{1,2}(\Qexn)^{*})} & \leq C\,.
\end{align*}
\end{lem}

The proof is given in Section \ref{sec:Proof-of-Lemma} (\thmref{Solution-Perforated-Domain}).

\medskip{}
The next step is passing to the limit $\eps\to0$. We do so in the
case that $\x$ is the realization of a stationary ergodic point process
$\X$ as defined below:
\begin{defn}[Random measure and shift-operator $\tau_{x}$]
\label{def:Random-Measure}\ \\
A \emph{random measure} $\mu_{\bullet}$ is a random variable with
values in $\M(\R^{d})$. It induces a probability distribution $\P$
on $\M(\Rd)$. Given the continuous map 
\begin{equation}
\tau_{x}:\,\M(\Rd)\to\M(\Rd)\,,\qquad\tau_{x}\xi(A):=\xi(A+x)\,,\label{eq:def-tau}
\end{equation}
a random measure is \emph{stationary} iff $\P(F)=\P(\tau_{x}F)$ for
every $F\in\B(\M(\Rd))$ and every $x\in\Rd$. In line with the above
setting, a random \emph{point process} $\X$ is a random measure with
$\P(\S(\R^{d}))=1$ and one quickly verifies that $\X$ is stationary
iff for every $N\in\N$, $x\in\Rd$ and bounded open $A\subset\Rd$
it holds 
\[
\P\left(\x\in\S(\R^{d}):\,\x\of A=N\right)=\P\left(\x\in\S(\R^{d}):\,\x\of{A+x}=N\right).
\]
We call a stationary random measure $\mu_{\bullet}$ \emph{ergodic}
iff the $\sigma$-algebra of $\tau$-invariant sets is trivial under
its distribution $\P$.
\end{defn}

\begin{rem}[Compatibility of thinning with shifts]
\ \\
The thinning map $\th$ is compatible with the shift $\tau_{x}$,
i.e. on $\S(\R^{d})$
\[
\th\circ\tau_{x}=\tau_{x}\circ\th\,.
\]
 
\end{rem}

\begin{lem}[Homogenized PDE for minimally smooth domains]
\label{lem:General-Homogenized-Limit}\ \\
Let $\X$ be a stationary ergodic point process and $n\in\N$ fixed.
For almost every realization $\x$ of $\X$, we have under \assuref{PDE-Parameters}:
\\
For $\eps>0$, let $\uexn\in L^{2}(I;W^{1,2}(\Qexn))$ be a solution
to \Eqref{System-in-Eps}. For a subsequence, there exist $\tilde{u}^{\eps}\in L^{2}(I;W^{1,2}(Q))$
with $\tilde{u}^{\eps}|_{\Qexn}=u^{\eps}$ such that $\tilde{u}^{\eps}\to u_{n}$
strongly in $L^{2}(I;L^{2}(Q))$ for some $u_{n}\in L^{2}(I;W^{1,2}(Q))$
with generalized time derivative $\del_{t}u_{n}\in L^{2}(I;W^{1,2}(Q)^{*})$.
\\
This $u_{n}$ is a weak solution to
\begin{align}
C_{1,\P^{(n)}}\del_{t}u_{n}-\nabla\cdot\big(A(u_{n})\A^{(n)}\nabla u_{n}\big)-C_{2,\P^{(n)}}h(u_{n}) & =C_{1,\P^{(n)}}f &  & \text{in }I\times Q\nonumber \\
A(u_{n})\A^{(n)}\nabla u_{n}\cdot\nu & =0 &  & \text{on }I\times\del Q\label{eq:Homogenized-System}\\
u_{n}(0,x) & =C_{1,\P^{(n)}}u_{0}(x) &  & \text{in }Q\nonumber 
\end{align}
with constants $C_{i,\P^{(n)}}>0$ depending on the distribution $\P^{(n)}$
of $\X^{(n)}$ and $\A^{(n)}$ being a symmetric positive semi-definite
matrix -- the so called effective conductivity based on the event
that the origin is not covered by $\bm\X^{(n)}$ (see \defref{Effective-Conductivity}).
\end{lem}

This is shown in Section \ref{sec:Proof-of-Lemma} (\thmref{Homogenized-System-for-n})
using two-scale convergence.

\subsection{Homogenization for irregular perforations}

When it comes to the final homogenization result, we will need the
following assumptions on the point process $\X$.
\begin{defn}[Admissible point process]
\label{def:Admissible-PP} \ \\
We call a point cloud $\x\in\S(\R^{d})$ \emph{admissible} iff the
following holds (with $r>0$ from \defref{boolean-model-filled}):
\begin{enumerate}
\item Equidistance Property: $\forall x,y\in\x:\,\left|x-y\right|\neq2r$.
\item Finite Clusters: For every $x\in\x$, we have that $\#\cC_{\x}(x)<\infty$.
\end{enumerate}
A stationary ergodic boundedly finite point process $\X$ is called
\emph{admissible} if its realizations are almost surely admissible\emph{. }
\end{defn}

\begin{defn}[Statistical connectedness]
\label{def:Statistically-Connected} \ \\
The random set $\bm\X^{\complement}$ is \emph{statistically connected}
iff the effective conductivity $\A$ (\defref{Effective-Conductivity})
based on the event that the origin is covered by $\bm\X^{\complement}$
is strictly positive definite.
\end{defn}

\begin{rem}[Sufficient condition for statistical connectedness]
\ \\
A criterion for statistical connectedness is given in Section \ref{sec:Sufficient-Condition-for-Effective-Conductivity},
namely the existence of sufficiently many so called percolation channels.
It also turns out that $\bm\X^{\complement}$ is statistically connected
if and only if the same holds for $\Xi\X^{\complement}$.
\end{rem}

\medskip{}
We may now state the main theorem of this work.
\begin{thm}[Homogenized limit for admissible point processes]
\label{thm:Main-Hom-Thm}~\\
Let $\X$ be an admissible point process and $\bm\X{}^{\complement}$
statistically connected. Under \assuref{PDE-Parameters}, we have
for almost every realization $\x$ of $\X$:\\
For every $n\in\N$, let $u_{n}$ be a homogenized limit from \lemref{General-Homogenized-Limit}.
Then, there exists a $u\in L^{2}(I;W^{1,2}(Q))$ with generalized
time-derivative $\del_{t}u\in L^{2}(I;W^{1,2}(Q)^{*})$ such that
for a subsequence
\begin{align*}
u_{n} & \xrightharpoonup[n\to\infty]{L^{2}(I;W^{1,2}(Q))}u\qquad\text{and}\qquad\del_{t}u_{n}\xrightharpoonup[n\to\infty]{L^{2}(I;W^{1,2}(Q)^{*})}\del_{t}u
\end{align*}
and $u$ is a weak solution to
\begin{align*}
C_{1,\P}\del_{t}u-\nabla\cdot\big(A(u)\A\nabla u\big)-C_{2,\P}h(u) & =C_{1,\P}f &  & \text{in }I\times Q\\
A(u)\A\nabla u\cdot\nu & =0 &  & \text{on }I\times\del Q\\
u(0,x) & =C_{1,\P}u_{0}(x) &  & \text{in }Q
\end{align*}
with constants $C_{i,\P}>0$ only depending on the distribution $\P$
of $\X$ and $\A$ being a symmetric positive definite matrix --
the so called effective conductivity $\A$ based on the event that
the origin is not covered by $\bm\X$ (\defref{Effective-Conductivity}).
In particular, the system does not depend on the chosen thinning procedure.
\end{thm}

\begin{proof}
This main theorem is proven in \thmref{Main-Theorem}.
\end{proof}
\begin{rem}[Random radii]
\ \\
Out of convenience, we have chosen the Boolean model with fixed radius
$r$ as our underlying model. One can easily generalize the procedure
to random independent radii $r\leq r_{\max}$. 
\end{rem}

\begin{rem}[Homogenization procedure]
\ \\
For fixed $\eps>0$, solutions $u^{\eps}=u_{\x}^{\eps}$ to \Eqref{System-in-Eps}
exist for admissible $\x\in\S(\R^{d})$ as $Q_{\xn}^{\eps}=Q_{\x}^{\eps}$
for $n$ large enough (\thmref{Approximation-Properties}). If $\x$
is a realization of some admissible point process $\X$, then this
is still not sufficient to pass to the limit $\eps\to0$. The missing
regularity of $\bm\X$ still prevents us from establishing a priori
estimates. \\
All in all, our procedure yields the following diagram:

{\Large{}
\[
\xymatrix{\xyR{3pc}\xyC{7pc}u^{\eps}\ar@{-->}[r]_{\eps\to\infty}^{?} & u\\
u_{n}^{\eps}\ar[u]^{n\uparrow\infty}\ar[r]_{\eps\to\infty} & u_{n}\ar[u]_{n\uparrow\infty}
}
\]
}Statistical connectedness of $\bm\X{}^{\complement}$ is crucial
to establish $W^{1,2}(Q)$-estimates for $u_{n}$. This indicates
that the direct limit passing $u^{\eps}\to u$ might only rely on
the statistical connectedness, but we cannot answer that as of yet.
\end{rem}

\subsection{Example: Poisson point processes}

In order to demonstrate that the class of point process satisfying
our assumptions is not empty, we show in Section \ref{sec:PPP} that
the Poisson point process $\X_{\poi}$ is indeed suitable for our
framework. We obtain the following. 
\begin{thm}[Admissibility and statistical connectedness for $\X_{\poi}$]
\ \\
In the subcritical regime (see \assuref{Subcritical-Regime}), we
have for the Poisson point process $\X_{\poi}$ that
\begin{itemize}
\item $\X_{\poi}$ is an admissible point process.
\item $\bm\X_{\poi}^{\complement}$ is statistically connected.
\end{itemize}
\end{thm}

While admissibility is easily proven, the statistical connectedness
is much harder to deal with. Most of Section \ref{sec:PPP} is dedicated
to this proof. It also builds up on Section \ref{sec:Sufficient-Condition-for-Effective-Conductivity}
in which we show that so called percolation channels yield statistical
connectedness.

\section{Thinning properties, surface measure and convergence of intensities
\label{sec:Properties-of-Thinning-Map}}

We first establish some properties of $\th:\ \S(\R^{d})\to\S(\R^{d})$,
most importantly the minimal smoothness of $\bm\x^{(n)}$.

\medskip{}

\begin{defn}[Minimal smoothness \cite{stein2016singular}]
\label{def:minimal-smooth}\ \\
An open set $P\subset\R^{d}$ is called \emph{minimally smooth} with
constants $(\delta,N,M)$ if we may cover $\del P$ by a countable
sequence of open sets $(U_{i})_{i}$ such that
\begin{enumerate}
\item $\forall x\in\R^{d}:\ \#\{U_{i}\,\vert\,x\in U_{i}\}\leq N$.
\item $\forall x\in\del P\,\exists U_{i}:\ \Ball{\delta}x\subset U_{i}$.
\item For every $i$, $\del P\cap U_{i}$ agrees (in some Cartesian system
of coordinates) with the graph of a Lipschitz function whose Lipschitz
semi-norm is at most $M$.
\end{enumerate}
\end{defn}

\begin{lem}[Uniform $\delta$ on individual clusters]
\label{lem:Admissible-PP-Property-3}\ \\
Let $\x\in\S(\R^{d})$ be an admissible point cloud. Then, for every
$x\in\x$
\[
\delta\of{\Xi\of{\cC_{\x}(x)}}>0\,.
\]
\end{lem}

\begin{proof}
Let $\x\in\S(\R^{d})$ and assume $\delta\of{\Xi\of{\cC_{\x}(x)}}=0$
for some $x\in\x$. Then, there must be some $p\in\del\Xi_{\x}$ with
$\delta(p)=0$. This together with bounded finiteness gives $x_{p},y_{p}\in\x$
such that $p\in B_{r}(x_{p})\cap B_{r}(y_{p})$, in particular $\biga x_{p}-y_{p}\biga=2r$.
This contradicts the equidistance property of $\x$.
\end{proof}
The thinning maps $\th$ have been constructed just to yield the following
theorem:
\begin{thm}[Minimal smoothness of thinned point clouds]
\label{thm:Minimal-Smoothness-By-F-n}\ \\
For every $\x\in\S(\R^{d})$, both $(\Xi\x^{(n)})^{\complement}$
and $(\bm\x^{(n)})^{\complement}$ are minimally smooth with $\delta=\frac{1}{n}$,
$M=\sqrt{2nr}$. Furthermore, every connected component of $\Xi\x^{(n)}$
or $\bm\x^{(n)}$ has diameter less than $2nr$.
\end{thm}

\begin{proof}
It remains to verify the estimate on $M$. Let $x=o=0_{\R^{d}}$ and
$y=\left(2r-n^{-1},0\dots,0\right)$. Then the Lipschitz constant
at the intersection of the two balls $\Ball rx$ and $\Ball ry$ is
less than $\sqrt{2nr}$.
\end{proof}
\begin{thm}[Further properties of $\th$]
\label{thm:Fn-compact-and-nice}~\\
The set $\S_{\A}(\Rd)$ of admissible point clouds is measurable in
the vague $\sigma$-algebra. Given $n\in\N$, it holds that $\th:\,\S(\Rd)\to\S(\Rd)$
is measurable, $\S^{(n)}:=\th\S(\Rd)$ is compact in the vague topology
and the following three properties of $\x\in\S(\Rd)$ are equivalent:
\begin{enumerate}
\item $\th\x=\x$ 
\item $\x\in\S^{(n)}$ 
\item (\ref{eq:thm:Fn-compact-and-nice-1})--(\ref{eq:thm:Fn-compact-and-nice-2})
hold:
\begin{align}
\forall x,y\in\x,x\neq y:\quad & d(x,y)\notin(0,\frac{1}{n})\cup(2r-\frac{1}{n},\,2r+\frac{1}{n})\,,\label{eq:thm:Fn-compact-and-nice-1}\\
\forall x\in\x:\quad & \#\mathcal{C}_{\x}(x)\leq n,\,\delta\of{\Ball r{\cC_{\x}(x)}}\geq\frac{1}{n}\,.\label{eq:thm:Fn-compact-and-nice-2}
\end{align}
 
\end{enumerate}
\end{thm}

\begin{proof}
$\th\x=\x$ implies $\x\in\S^{(n)}$ since $\th\x\in\S^{(n)}$ and
vice versa $\x\in\S^{(n)}$ implies $\th\x=\x$ by definition of $\th$.
By construction of $\th$ it follows that (\ref{eq:thm:Fn-compact-and-nice-1})--(\ref{eq:thm:Fn-compact-and-nice-2})
hold if and only if $\x\in\S^{(n)}$. 

Consider the space of (non-simple) counting measures $\cN(\R^{d})\subset\M(\R^{d})$,
i.e.
\[
\cN(\R^{d}):=\big\{\mu\in\M(\R^{d})\,\vert\,\text{\ensuremath{\mu=\sum_{k\in\mathcal{I}\subset\N}a_{k}\delta_{x_{k}}\text{ such that }}}a_{k}\in\N\text{ and }x_{k}\in\R^{d}\big\}\,.
\]
We see, e.g. in \cite{MR2371524}, that 
\begin{itemize}
\item $\S(\R^{d})$ and $\cN(\R^{d})$ are both measurable w.r.t. the Borel-$\sigma$-algebra
of $\M(\R^{d})$.
\item $\S(\R^{d})\subset\cN(\R^{d})$ and $\cN(\R^{d})$ is closed in $\M(\R^{d})$.
In particular, $\cN(\R^{d})$ is also complete under the Prohorov
metric.
\end{itemize}
Now $\S^{(n)}$ is precompact because of the characterization of precompact
sets in the vague topology: For every bounded open $A\subset\Rd$,
it holds that $\sup_{\x\in\S^{(n)}}\x(A)\leq C\left(\mathrm{diam}\,A\right)^{d}$
with $C$ depending only on $n$. It remains to show that $\S^{(n)}$
is closed as a subset of $\cN(\R^{d})$. Let $\left(\x_{j}\right)_{j\in\N}\subset\S^{(n)}$
be a converging sequence with limit $\x\in\cN(\Rd)$. One checks that
\Eqref{thm:Fn-compact-and-nice-1} (namely $d(x,y)\notin\big(0,\frac{1}{n}\big)$)
ensures $\x\in\S(\R^{d})$, e.g. in a procedure similar to the proof
of \cite[Lemma 9.1.V]{MR2371524}. We observe that for every $x,y\in\x$,
there exist $x_{j},y_{j}\in\x_{j}$ such that $x_{j}\to x$, $y_{j}\to y$
as $j\to\infty$. This implies by a limit in (\ref{eq:thm:Fn-compact-and-nice-1})
that $\x$ still satisfies (\ref{eq:thm:Fn-compact-and-nice-1}).

For $x\in\x$, one checks that \Eqref{thm:Fn-compact-and-nice-1}
(namely $d(x,y)\notin\big(2r-\frac{1}{n},2r+\frac{1}{n}\big)$) implies
$\#\cC_{\x}(x)\leq n$.

Let $p\in\partial\Xi(\x)$ and let $\left\{ x^{(1)},\dots x^{(K)}\right\} =\Ball{10r}p\cap\x$
with sequences $x_{j}^{(k)}\to x^{(k)}$, $x_{j}^{(k)}\in\x_{j}$.
Given $\eta>0$, let $J\in\N$ such that for all $j>J$ and $k=1,\dots K$
it holds $\left|x^{(k)}-x_{j}^{(k)}\right|<\eta$. Then there exists
$p_{j}\in\partial\Xi(\x_{j})$ such that $\left|p_{j}-p\right|<\eta$
and $\partial\Xi(\x_{j})$ is a Lipschitz graph in the ball $\Ball{2\delta}{p_{j}}$
for every $\delta<\frac{1}{n}$. Hence $\partial\Xi(\x_{j})$ is a
Lipschitz graph in the ball $\Ball{2\delta-\eta}p$. Because the Lipschitz
regularity of $\partial\Xi(\x_{j})$ changes continuously under slight
shifts of the balls, there exists $\eta_{0}$ such that for $\eta<\eta_{0}$
and $\partial\Xi(\x)$ is Lipschitz graph in $\Ball{2\delta-2\eta}p$.
Since $\eta$ is arbitrary, we find $\partial\Xi(\x)$ is Lipschitz
graph in $\Ball{2\delta}p$ for every $\delta<n^{-1}$, implying $\delta(p)\geq\frac{1}{n}$.
Since this holds for every $p$, we conclude (\ref{eq:thm:Fn-compact-and-nice-2})
and $\S^{(n)}$ is compact.

To see that $\S_{\A}(\Rd)$ is measurable, consider for $\x\in\S(\Rd)$
\begin{align*}
\x_{1,m,R} & :=\left\{ x\in\x\,\vert\,x\notin\Ball Ro\text{ or }\,d(x,y)\notin(0,\frac{1}{m})\cup(2r-\frac{1}{m},\,2r+\frac{1}{m})\ \forall y\in\x\right\} \,,\\
\x_{2,l,R} & :=\left\{ x\in\x\,\vert\,x\notin\overline{\Ball Ro}\text{ or }\#\mathcal{C}_{\x}(x)\leq l,\,\delta\of{\Ball r{\cC_{\x}(x)}}>\frac{1}{l}\right\} \,,
\end{align*}
and define
\begin{align*}
\Fth_{1,m,R}\x & :=\x_{1,m,R} & \S^{(1,m,R)} & :=\Fth_{1,m,R}\S(\R^{d})\\
\Fth_{2,l,R}\x & :=\x_{2,l,R}\,. & \S^{(2,m,R)} & :=\Fth_{2,m,R}\S(\R^{d})\,.
\end{align*}
We check that $\S^{(1,m,R)}$ is a closed subset inside $\S(\R^{d})$
(repeat the arguments above), i.e. $\overline{\S^{(1,m,R)}}\cap\S(\R^{d})=\S^{(1,m,R)}$.
In particular, $\S^{(1,m,R)}$ is measurable w.r.t. the vague topology
of $\M(\R^{d})$. Similarly, one shows that $\S(\R^{d})\backslash\S^{(2,m,R)}$
is closed as a subset inside $\S(\R^{d})$. Again, this shows that
$\S^{(2,m,R)}$ is measurable. Consider now the measurable sets
\[
\S^{(1,\infty,\infty)}:=\bigcap_{R\in\N}\bigcup_{m\in\N}\S^{(1,m,R)}\qquad\text{and}\qquad\S^{(2,\infty,\infty)}:=\bigcap_{R\in\N}\bigcup_{m\in\N}\S^{(2,m,R)}\,.
\]
We see that 
\begin{enumerate}
\item $\x\in\S^{(1,\infty,\infty)}$ if and only if for all $x,y\in\x$
, it holds $d(x,y)\neq r$.
\item $\x\in\S^{(2,\infty,\infty)}$ if and only if for every $x\in\x$,
it holds that $\#\mathcal{C}_{\x}(x)<\infty$ and $\delta\of{\Ball r{\cC_{\x}(x)}}>0$. 
\end{enumerate}
Therefore,
\[
\S_{\A}(\Rd)=\S^{(1,\infty,\infty)}\cap\S^{(2,\infty,\infty)}
\]
is measurable. 

To see that $\th:\ \S(\Rd)\to\S(\Rd)$ is measurable, recall $\th=\Fth_{2,n}\circ\Fth_{1,n}$
from \defref{Thinning-Maps}. It therefore suffices to show that the
following maps are measurable:
\[
\Fth_{1,n}:\ \S(\R^{d})\to\Fth_{1,n}\S(\R^{d})\subset\S(\R^{d})\qquad\text{and}\qquad\Fth_{2,n}:\ \Fth_{1,n}\S(\R^{d})\to\S(\R^{d})\,.
\]
For $f\in C_{c}(\Rd)$, consider the evaluation by $f$, i.e. 
\[
M_{f}:\,\S(\Rd)\to\R\,,\qquad\x\mapsto\int_{\Rd}f\,\d\x\,.
\]
If $f\geq0$, we observe the upper semi-continuity of
\[
M_{f}\circ\Fth_{1,n}:\ \S(\R^{d})\to\R\qquad\text{and}\qquad M_{f}\circ\Fth_{2,n}:\ \Fth_{1,n}\S(\R^{d})\to\R\,.
\]
We have lower semi continuity for $f\leq0$ since $M_{-f}=-M_{f}$.
 Therefore $M_{f}\circ\Fth_{i,n}$ with $i\in\{1,2\}$ is measurable
in the cases $f\geq0$ and $f\leq0$ and hence in general. Since the
vague topology is generated by $\big(M_{f}\big)_{f\in C_{c}(\R^{d})}$,
we conclude that $\Fth_{1,n}$ and $\Fth_{2,n}$ are measurable.
\end{proof}

\begin{rem}[Fine details of \thmref{Fn-compact-and-nice}]
\ 
\begin{itemize}
\item For $\S^{(n)}:=\th(\S(\R^{d})$), we have that
\[
\bigcup_{n\in\N}\S^{(n)}\subsetneq\S_{\A}(\R^{d})\subsetneq\big\{\x\,\vert\,\lim_{n\to\infty}\th\x=\x\big\}\subsetneq\S(\R^{d})\subsetneq\overline{\bigcup_{n\in\N}\S^{(n)}}=\cN(\R^{d})\,.
\]
\item $M_{f}\circ\Fth_{2,n}$ is \emph{not} upper semi-continuous on $\S(\R^{d})$
(in contrast to $\Fth_{1,n}\S(\R^{d})$): The condition that $d(x,y)\notin(2r-\frac{1}{n},\,2r+\frac{1}{n})\ \forall x,y\in\Fth_{1,n}\x$
is crucial to ensure that clusters do not change sizes.
\end{itemize}
\end{rem}

\begin{defn}
\label{def:G_n}We define the events that the origin is not covered
by the filled-up Boolean model, i.e.
\[
\bG:=\left\{ \x\in\S(\R^{d})\,\vert\,o\notin\bm\x\right\} \qquad\text{and}\qquad\bG_{n}:=\th{}^{-1}(\bG)=\left\{ \x\in\S(\R^{d})\,\vert\,o\notin\bm\x^{(n)}\right\} \,.
\]
This gives us for $x\in\R^{d}$ that
\[
\I_{\bm\x^{\complement}}(x)=\I_{\bG}(\tau_{x}\x)\,.
\]
\medskip{}
We will later consider the effective conductivities based on these
events.
\end{defn}

\begin{thm}[Approximation properties]
\label{thm:Approximation-Properties}\ \\
Let $\x\in\S(\R^{d})$ be an admissible point cloud.
\begin{enumerate}
\item For every bounded domain $\Lambda$, there exists an $N(\x,\Lambda)\in\N$
such that for every $n\geq N(\x,\Lambda)$
\[
\x^{(n)}\cap\Lambda=\x\cap\Lambda\,,\qquad\text{in particular}\qquad\x=\bigcup_{n\in\N}\x^{(n)}\,.
\]
\item For every bounded domain $\Lambda$, there exists an $\tilde{N}(\x,\Lambda)\in\N$
such that for every $n\geq\tilde{N}(\x,\Lambda)$
\[
\bm\x{}^{(n)}\cap\Lambda=\bm\x\cap\Lambda,\qquad\text{in particular}\qquad\bm\x=\bigcup_{n\in\N}\bm\x{}^{(n)}.
\]
\item There exists an $N=N(\x)\in\N$ such that for every $n\geq N$:
\[
o\notin\bm\x^{(n)}\iff o\notin\bm\x\,.
\]
In particular, $\bigcap_{n\in\N}\bG_{n}\backslash\bG$ only consists
of non-admissible point clouds.
\end{enumerate}
\end{thm}

\begin{proof}
\ 
\begin{enumerate}
\item Boundedness of $\Lambda$ implies that there are only finitely many
mutually disjoint clusters $\cC_{\x}(x_{i})$, $i=1,\dots,N_{\cC}$
that intersect with $\Lambda$. Furthermore, because $\#\left(\x\cap\Ball r{\Lambda}\right)<\infty$
and because of Property 1 of admissible point clouds, we know 
\[
\min\left\{ \left|\left|x-y\right|-2r\right|:\,x,y\in\x\cap\Ball r{\Lambda},\,x\neq y\right\} >0
\]
and \lemref{Admissible-PP-Property-3} yields
\[
\min\left\{ \delta\of{\Xi\cC_{\x}(x_{i})}:\,\cC_{\x}(x_{i})\cap\Lambda\neq\emptyset\right\} >0\,.
\]
This implies the first statement.
\item By making $\Lambda$ larger, we may assume $\Lambda=[-k,k]^{d}$ for
some $k\in\N$. For $n\geq N(\Ball r{[-k,k]^{d}})$: 
\[
[-k,k]^{d}\backslash\Xi\x^{(n)}=[-k,k]^{d}\backslash\Xi\x.
\]
$[-k,k]^{d}\backslash\Xi\x$ only has finitely many connected components
$\mathcal{C}_{i}$. Take one of these connected components $\mathcal{C}_{i}$
and suppose it lies in $\bm\x$. Then, it has to be encircled by finitely
many balls $\Ball rx$ in $\Xi\x$. Let $n_{i}$ large enough such
that all these $x$ lie in $\x^{(n_{i})}$. Then, $\mathcal{C}_{i}\subset\bm\x^{(n_{i})}$.
We may do so for every $\mathcal{C}_{i}$. Take 
\[
\tilde{N}(\x,\Lambda):=\max\big\{ n_{i},\,N(\Ball r{[-k,k]^{d}})\big\}.
\]
For every $n\geq\tilde{N}(\x,\Lambda)$, the connected components
$\mathcal{C}_{i}$ of $[-k,k]^{d}\backslash\Xi\x^{(n)}$ and $[-k,k]^{d}\backslash\Xi\x$
are identical since $[-k,k]^{d}\backslash\Xi\x^{(n)}=[-k,k]^{d}\backslash\Xi\x$.
Therefore, we get the claim
\begin{align*}
[-k,k]^{d}\backslash\bm\x^{(n)} & =\Big([-k,k]^{d}\backslash\Xi\x^{(n)}\Big)\backslash\bigcup_{C_{i}\subset\bm\x^{(n)}}\mathcal{C}_{i}\\
 & =\Big([-k,k]^{d}\backslash\Xi\x\Big)\backslash\bigcup_{C_{i}\subset\bm\x}\mathcal{C}_{i}\quad=[-k,k]^{d}\backslash\bm\x\,.
\end{align*}
\item This is a direct consequence of Point 2. If $\x\in\bigcap_{n\in\N}\bG_{n}\backslash\bG$,
then $o\notin\bm\x^{(n)}$ for every $n$ but $o\in\bm\x$. Therefore,
$\x$ cannot be admissible by Point 2.
\end{enumerate}
\end{proof}
\begin{defn}[Surface measure of $\bm\x$]
\label{def:Surface-Measure-of-BM}\ \\
We define the \emph{surface measure} for $\x\in\S(\R^{d})$
\[
\mu_{\x}(A):=\H_{\llcorner\del\bm\x}^{d-1}(A)=\H^{d-1}(A\cap\del\bm\x)\,.
\]
Note that $\mu_{\x}\of{[0,1]^{d}}\leq\H^{d-1}\of{\Ball ro}\,\cdot\,\x\of{\Ball r{[0,1]^{d}}}$.
\end{defn}

\begin{rem}[Distributions of $\X^{(n)},\,\mu_{\bullet}$]
\ \\
Given a point process $\X$, we can consider $\X^{(n)}$ and $\mu_{\X}$.
Both come with their own distributions, but they are still driven
by $\X$ in a $\tau$-compatible way. Therefore, we can express their
distributions and all relevant quantities in terms of the distribution
$\P$ of $\X$. For example, the distribution of $\X^{(n)}$ is $\P\circ\th{}^{-1}$.\\
\end{rem}

\begin{defn}[Intensity of random measure]
\ \\
Given a stationary random measure $\tilde{\mu}$, we define its intensity
to be
\[
\lambda(\tilde{\mu}):=\E\big[\tilde{\mu}\big([0,1)^{d}\big)\big]\,.
\]
\end{defn}

\begin{lem}[Convergence of intensities]
\label{lem:Convergence-of-Intensities}\ \\
Let $\X$ be an admissible point process with finite intensity $\lambda(\X)$.
Then,
\[
\lim_{n\to\infty}\lambda(\X^{(n)})=\lambda(\X)\qquad\text{and}\qquad\lim_{n\to\infty}\lambda(\mu_{\X^{(n)}})=\lambda(\mu_{\X})\,.
\]
\end{lem}

\begin{proof}
``Almost surely'' is to be understood w.r.t. the distribution $\P$
of $\X$.
\begin{enumerate}
\item By \thmref{Approximation-Properties}, we have almost surely $\X^{(n)}\of{[0,1]^{d}}\to\X\of{[0,1]^{d}}$
as $n\to\infty$. Dominated convergence with majorant $\X\of{[0,1]^{d}}$
yields
\[
\lambda(\X^{(n)})=\E[\x^{(n)}\of{[0,1]^{d}}]\to\E[\x\of{[0,1]^{d}}]=\lambda(\X).
\]
\item Again, by \thmref{Approximation-Properties}, we have almost surely
$\bm\X^{(n)}\cap[0,1]^{d}\to\bm\X\cap[0,1]^{d}$, in particular
\[
\mu_{\X^{(n)}}\of{[0,1]^{d}}=\H_{\llcorner\del\bm\X^{(n)}}^{d-1}\of{[0,1]^{d}}\to\H_{\llcorner\del\bm\X}^{d-1}\of{[0,1]^{d}}=\mu_{\X}\of{[0,1]^{d}}.
\]
Dominated convergence yields again convergence of intensities.
\end{enumerate}
\end{proof}
\begin{rem}[Local convergence]
\ \\
The convergence in \thmref{Approximation-Properties} is much stronger
than what is actually needed to prove the convergence of intensities.
Indeed, we could prove convergence even for so called tame and local
functions $f:\ \M(\R^{d})\to\R$ for which the intensity $\lambda$
is just one special case $f(\x):=\x\of{[0,1]^{d}}$.
\end{rem}

\section{Effective conductivity and cell solutions \label{sec:Convergences-of-Cell-Solutions}}

The structure $(\M(\R^{d}),\B(\M(\R^{d})),\P,\tau)$ as in \defref{Random-Measure}
is a so called dynamical system:
\begin{defn}[Dynamical system, stationarity, ergodicity]
\label{def:Dynamic-System}\ \\
Let $(\Omega,\F,\p)$ be a separable metric probability space. A \emph{dynamical
system} $\tau=(\tau_{x})_{x\in\R^{d}}$ is a family of measurable
mappings $\tau_{x}:\ \Omega\to\Omega$ satisfying 
\begin{itemize}
\item Group property: \\
$\tau_{0}=\text{id}_{\Omega}$ and $\tau_{x+y}=\tau_{x}\circ\tau_{y}$
for any $x,\,y\in\mathbb{R}^{d}$.
\item Measure preserving:\\
For any $x\in\mathbb{R}^{d}$ and any $F\in\mathcal{F}$, we have
$\p\left(\tau_{x}(F)\right)=\p\left(F\right)$.
\item Continuity:\\
The map $\mathcal{T}:\ \Omega\times\mathbb{R}^{d}\to\Omega$, $\left(\omega,\,x\right)\mapsto\tau_{x}(\omega)$
is continuous w.r.t. the product topology on $\Omega\times\R^{d}$.
\end{itemize}
$\tau$ is called \emph{ergodic} if the $\sigma$-algebra of $\tau$-invariant
sets is trivial under $\p$.
\end{defn}

\medskip{}
Our practical setting will always be some $\Omega\subset\M(\R^{d})$,
but we will still work with abstract dynamical systems in Section
\ref{sec:Convergences-of-Cell-Solutions} and \ref{sec:Proof-of-Lemma}.

\subsection{Potentials and solenoidals}

Let $(\Omega,\B(\Omega),\p,\tau)$ be a dynamical system. We write
$L^{2}(\Omega):=L^{2}(\Omega,\,\p)$. The dynamical system $\tau$
introduces a strongly continuous group action on $L^{2}(\Omega)\to L^{2}(\Omega)$
through $T_{x}f(\omega):=f(\tau_{x}\omega)$ with the $d$ independent
generators 
\[
\rmD_{i}f:=\lim_{t\to0}\frac{1}{t}\left(f-f(\tau_{te_{i}}\bullet)\right)
\]
with domain $\mathcal{D}_{i}$ where $\left(e_{i}\right)_{i=1,\dots,d}\subset\R^{d}$
is the canonical Euclidean basis. Introducing 
\[
H^{1}(\Omega):=\bigcap_{i=1}^{d}\mathcal{D}_{i}\subset L^{2}(\Omega)
\]

\noindent and the gradient $\nabla_{\omega}f:=\left(\rmD_{1}f,\dots,\rmD_{d}f\right)^{\top}$,
we can define the space of potential vector fields
\[
\mathcal{V}_{\pot}^{2}(\Omega):=\big\{\nabla_{\tilde{\omega}}f\,\vert\,f\in H^{1}(\Omega)\text{ and }\int_{\Omega}\nabla_{\tilde{\omega}}f\dx\p(\omega)=0_{\R^{d}}\big\}\,.
\]
Defining $L_{\sol}^{2}(\Omega):=\cV_{\pot}^{2}(\Omega)^{\bot}$, we
find with $u_{\omega}(x):=u(\tau_{x}\omega)$
\begin{align}
L_{\pot}^{2}(\Omega) & :=\left\{ u\in L^{2}(\Omega;\Rd)\,:\,u_{\omega}\in L_{\pot,\loc}^{2}(\Rd)\,\,\mbox{for }\p-\mbox{a.e. }\omega\in\Omega\right\} \,,\nonumber \\
L_{\sol}^{2}(\Omega) & =\left\{ u\in L^{2}(\Omega;\Rd)\,:\,u_{\omega}\in L_{\sol,\loc}^{2}(\Rd)\,\,\mbox{for }\p-\mbox{a.e. }\omega\in\Omega\right\} \,,\label{eq:sto-Lp-pot-sol-omega-general-1}\\
\cV_{\pot}^{2}(\Omega) & =\big\{ u\in L_{\pot}^{2}(\Omega)\,:\,\int_{\Omega}u\,\d\p=0\big\}\,,\nonumber 
\end{align}
because $\Omega$ is separable metric and where
\begin{align*}
L_{\pot,\loc}^{2}(\Rd) & :=\left\{ u\in L_{\loc}^{2}(\Rd;\Rd)\,\,|\,\,\forall U\,\,\mbox{bounded domain }\exists\varphi\in W^{1,2}(U)\,:\,u=\nabla\varphi\right\} \,,\\
L_{\sol,\loc}^{2}(\Rd) & :=\big\{ u\in L_{\loc}^{2}(\Rd;\Rd)\,\,|\,\,\int_{\Rd}u\cdot\nabla\varphi=0\,\,\forall\varphi\in C_{c}^{1}(\Rd)\big\}\,.
\end{align*}
For $A\subset\Omega$ measurable, we define
\[
\cV_{\pot}^{2}(A\vert\Omega):=cl_{L^{2}(A)^{d}}\left\{ v_{\vert_{A}}\,\vert\,v\in\mathcal{V}_{\pot}^{2}(\Omega)\right\} .
\]

\subsection{Cell solutions and effective conductivity}
\begin{defn}[Cell solutions]
\label{def:Cell-Solutions}\ \\
Let $(\Omega,\B(\Omega),\p)$ be a separable metric probability space
with dynamical system $\tau$ and let $\O\in\B(\Omega)$. The $i$-th
\emph{cell solution} $w_{i}\in\cV_{\pot}^{2}(\O\vert\Omega)$ is the
unique solution (after the Riesz representation theorem) of
\[
\forall v\in\mathcal{V}_{\pot}^{2}(\Omega):\qquad\int_{\O}[w_{i}+e_{i}]\cdot v\dx\p(\omega)=0\,.
\]
The cell solutions satisfy
\[
\bign w_{i}\bign_{L^{2}(\O)^{d}}\leq\sqrt{\p(\O)}\leq1\,,
\]
and can be grouped in the matrix 
\[
W_{\O}:=(w_{1},\dots,w_{d})\,.
\]
\end{defn}

\begin{defn}[Effective conductivity $\A$]
\label{def:Effective-Conductivity}\ \\
Let $w_{i}$ the cell solution on $\O\in\B(\Omega)$. The \emph{effective
conductivity} $\mathcal{A}$ based on the event $\O$ is defined as
\begin{align}
\A & :=\int_{\O}(I_{d}+W_{\O})^{t}(I_{d}+W_{\O})\dx\p(\omega)\,.\label{eq:eff-cond-infty}
\end{align}
with $I_{d}$ being the identity matrix. We observe for the entries
$\left(\mathcal{A}_{i,j}\right)_{i,j=1,\dots,d}$ of $\A$ that 
\begin{align}
\mathcal{A}_{i,j} & =\int_{\O}\left[e_{i}+w_{i}(\omega)\right]\cdot\left[e_{j}+w_{j}(\omega)\right]\dx\p(\omega)=\int_{\O}\left[e_{i}+w_{i}(\omega)\right]\cdot e_{j}\dx\p(\omega)\,.\label{eq:eff-cond-linear}
\end{align}
We write $\alpha_{\A}\geq0$ for its smallest eigenvalue.
\end{defn}

\begin{thm}[Convergence of cell solutions]
\ \\
Let $\left(\O_{n}\right)_{n\in\N}\subset\B(\Omega)$ and $\O\in\B(\Omega)$
such that $\I_{\O_{n}}\to\I_{\O}$ $\p$-almost surely and $\O_{n}\supset\O$
for every $n\in\N$. The sequence of cell solutions $w_{i}^{(n)}$
to the cell problem on $\O_{n}$ satisfies 
\[
w_{i}^{(n)}\weakto w_{i}\qquad\text{as }n\to\infty\,,
\]
where $w_{i}\in\cV_{\pot}^{2}(\O\vert\Omega)$ is the $i$-th cell
solution on $\O$
\[
\forall v\in\cV_{\pot}^{2}(\Omega):\qquad\int_{\O}\left[w_{i}+e_{i}\right]\cdot v\dx\p(\omega)=0\,.
\]
\end{thm}

\begin{proof}
\noindent We first check that the limit satisfies $\int_{\O}\left[w_{i}+e_{i}\right]\cdot v\dx\p(\omega)=0$
and then $w_{i}\in\cV_{\pot}^{2}(\O\vert\Omega)$.

\noindent 1. The a priori estimate yields a $L^{2}$-weakly convergent
subsequence of $w_{i}^{(n)}\rightharpoonup w_{i}\in L^{2}(\Omega)^{d}$
after extending $w_{i}^{(n)}$ to the whole of $\Omega$ via 0. Let
$v\in\cV_{\pot}^{2}(\Omega)$. We have $\I_{\O_{n}}\to\I_{\O}$ $\p$-almost
surely, so dominated convergence yields
\begin{align*}
\lim_{n\to\infty}\int_{\O_{n}}e_{i}\cdot v\dx\p(\omega) & =\int_{\O}e_{i}\cdot v\dx\p(\omega)
\end{align*}
while weak convergence yields
\begin{align*}
\lim_{n\to\infty}\int_{\O_{n}}w_{i}^{(n)}\cdot v\dx\p(\omega) & =\int_{\Omega}w_{i}\cdot v\dx\p(\omega)
\end{align*}
We also have
\[
\I_{\O_{n}}w_{i}^{(n)}=w_{i}^{(n)}\xrightharpoonup[n\to\infty]{L^{2}(\Omega)}w_{i}\,,
\]
which implies
\[
\I_{\O}w_{i}=w_{i}.
\]
Therefore, with $w_{i}\in L^{2}(\O)^{d}$:
\begin{align*}
0 & =\lim_{n\to\infty}\int_{\O_{n}}[w_{i}^{(n)}+e_{i}]\cdot v\dx\p(\omega)=\int_{\O}\left[w_{i}+e_{i}\right]\cdot v\dx\p(\omega)\,.
\end{align*}

2. $\cV_{\pot}^{2}(\O\vert\Omega)\subset L^{2}(\O)^{d}$ is closed
and convex, so it is also weakly closed. We construct a weakly converging
sequence in $\cV_{\pot}^{2}(\O\vert\Omega)$ that converges to $w_{i}$.
Since $w_{i}^{(n)}\in\cV_{\pot}^{2}(\O_{n}\vert\Omega)$, we find
$v^{(n)}\in\mathcal{V}_{\pot}^{2}(\Omega)$ such that
\[
\|w_{i}^{(n)}-\I_{\O_{n}}v^{(n)}\|_{L^{2}(\O)^{d}}\leq\frac{1}{n}\,.
\]
Since $w_{i}^{(n)}\rightharpoonup w_{i}$, we get
\[
\I_{\O_{n}}v^{(n)}\xrightharpoonup[n\to\infty]{L^{2}(\Omega)^{d}}w_{i}\,.
\]
Note that $\left(\I_{\O_{n}}-\I_{\O}\right)v^{(n)}$ is a bounded
sequence that is point-wise convergent to 0 because of the weak convergence
above and $\O_{n}\supseteq\O$. Therefore, it is weakly convergent
to 0 and we obtain
\[
\I_{\O}v^{(n)}\xrightharpoonup[n\to\infty]{L^{2}(\Omega)^{d}}w_{i}\,.
\]
$\I_{\O}v^{(n)}\in\cV_{\pot}^{2}(\O\vert\Omega)$, so we get that
$w_{i}\in\cV_{\pot}^{2}(\O\vert\Omega)$.
\end{proof}
\begin{cor}[Convergence of effective conductivities]
\label{cor:Convergence-of-Effective-Conductivities}\ \\
Let $\left(\O_{n}\right)_{n\in\N}\subset\B(\Omega)$ and $\O\in\B(\Omega)$
such that $\I_{\O_{n}}\to\I_{\O}$ $\p$-almost surely and $\O_{n}\supset\O$
for every $n\in\N$. Let $\A^{(n)}$ be the effective conductivity
of $\O_{n}$ and $\A$ be the effective conductivity of $\O$. Then,
\[
\mathcal{A}^{(n)}\xrightarrow{n\to\infty}\mathcal{A}\,.
\]
\end{cor}

\begin{proof}
This follows from \Eqref{eff-cond-linear} and weak convergence $w_{i}^{(n)}\rightharpoonup w_{i}$.
\end{proof}
\begin{rem}[Variational formulation]
\ \\
There is another way to define $\A$: For $\eta\in\R^{d}$, $W_{\O}\eta$
(see \defref{Cell-Solutions}) is the unique minimizer to
\[
\min_{\varphi\in\cV_{\pot}^{2}(\O)}\int_{\O}\biga\eta+\varphi\biga^{2}\dx\p(\omega)
\]
and therefore
\[
\eta^{t}\A\eta=\int_{\O}\biga(I_{d}+W_{\O})\eta\biga^{2}\dx\p(\omega)=\min_{\varphi\in\cV_{\pot}^{2}(\O)}\int_{\O}\biga\eta+\varphi\biga^{2}\dx\p(\omega).
\]
This equality is related to \thmref{Variational-Formulation-Effective-Conductivity}.
\end{rem}

\subsection{Pull-back for thinning maps \label{subsec:Pull-back}}

In Section \ref{sec:Proof-of-Lemma}, we will use two-scale convergence
to homogenize \Eqref{System-in-Eps} for fixed $n$. This process
is more convenient to handle if the underlying probability space is
compact. Here we show that we may take $\th\S(\R^{d})$ as the underlying
probability space instead of $\S(\R^{d})$.
\begin{thm}
\label{thm:pull-back-cell-solution}\ \\
Let $\p$ be a distribution on $\S(\Rd)$ and let $\S^{(n)}:=\th\S(\R^{d})$
with the push-forward measure $\tilde{\p}_{n}:=\p\circ\th^{-1}$.
Recall $\bG_{n}:=\left\{ \x\in\S(\R^{d}):\,o\notin\bm\x^{(n)}\right\} $
and let $\tilde{\bG}_{n}:=\left\{ \x\in\S^{(n)}:\,o\notin\bm\x\right\} =\th\bG_{n}$.
Let $w_{i}^{(n)}$ be the cell solutions on $\bG_{n}$ and $\tilde{w}_{i}^{(n)}$
the cell solutions on $\tilde{\bG}_{n}$ for their respective dynamical
systems. Then, for every $i,j\in\left\{ 1,\dots,d\right\} $, it holds
that
\begin{equation}
\int_{\bG_{n}}\big[w_{i}^{(n)}+e_{i}\big]\cdot e_{j}\dx\p=\int_{\tilde{\bG}_{n}}\big[\tilde{w}_{i}^{(n)}+e_{i}\big]\cdot e_{j}\dx\tilde{\p}_{n}\,.\label{eq:thm:pull-back-cell-solution}
\end{equation}
\end{thm}

\begin{lem}[Properties of pull-back functions]
\label{lem:pull-back} \ \\
Let $(\Omega,\F,\p,\tau)$, $(\tilde{\Omega},\tilde{\F},\tilde{\p},\tilde{\tau})$
be dynamical systems, $\phi:\ \Omega\to\tilde{\Omega}$ measurable
such that $\tilde{\p}=\p\circ\phi^{-1}$ and such that for every $x\in\R^{d}$
\begin{equation}
\phi\circ\tau_{x}=\tilde{\tau}_{x}\circ\phi\,.\label{eq:Pullback-Compatible}
\end{equation}
Then, the following holds:\\
For every $\tilde{f}\in L^{2}(\tilde{\Omega})^{d}$, we have $f:=\tilde{f}\circ\phi\in L^{2}(\Omega)^{d}$
with $\bign f\bign_{L^{2}(\Omega)^{d}}=\bign\tilde{f}\bign_{L^{2}(\tilde{\Omega})^{d}}$.
If $\tilde{f}\in\cV_{\pot}^{2}(\tilde{\Omega})$, then $f\in\cV_{\pot}^{2}(\Omega)$.
If $\tilde{f}\in L_{\sol}^{2}(\tilde{\Omega})$, then $f\in L_{\sol}^{2}(\Omega)$.
$f$ is called the \emph{pull-back }of $\tilde{f}$.
\end{lem}

\begin{proof}
Due to $\tilde{\p}=\p\circ\phi^{-1}$, we immediately obtain for arbitrary
measurable $\tilde{f}\in L^{1}(\tilde{\Omega})^{d}$ and its pull-back
$f$
\begin{align}
\int_{\tilde{\Omega}}\tilde{f}\dx\tilde{\p} & =\int_{\Omega}\tilde{f}\circ\phi\dx\p=\int_{\Omega}f\dx\p\,.\label{eq:Pull-Back-Integral-Equation}
\end{align}
Therefore $\bign f\bign_{L^{2}(\Omega)^{d}}=\bign\tilde{f}\bign_{L^{2}(\tilde{\Omega})^{d}}$
and \Eqref{Pullback-Compatible} yields $\tilde{f}\in\cV_{\pot}^{2}(\tilde{\Omega})\implies f\in\cV_{\pot}^{2}(\Omega)$.
For $\tilde{f}\in L_{\sol}^{2}(\tilde{\Omega})$, $f\in L_{\sol}^{2}(\Omega)$
follows from $\phi\circ\tau_{x}=\tilde{\tau}_{x}\circ\phi$ and checking

\[
\int_{\Lambda}f(\tau_{x}\omega)\cdot\nabla\varphi(x)\dx x=\int_{\Lambda}\tilde{f}(\tilde{\tau}_{x}\omega)\cdot\nabla\varphi(x)\dx x=0
\]
for $\p$-almost every $\omega$ and every $\varphi\in C_{c}^{1}(\Lambda)$
on a bounded domain $\Lambda\subset\R^{d}$.
\end{proof}
\begin{proof}[Proof of Theorem \ref{thm:pull-back-cell-solution}]
Let $\overline{w}_{i}$ be the pull-back of $\tilde{w}_{i}^{(n)}$
according to Lemma \ref{lem:pull-back}. We see $\th^{-1}\tilde{\bG}_{n}=\bG_{n}$,
so $\overline{w}_{i}$ has support in $\bG_{n}$. Let $\tilde{v}_{k}\in\cV_{\pot}^{2}(\S^{(n)})$
with $\bign\tilde{w}_{i}^{(n)}-\tilde{v}_{k}\bign_{L^{2}(\tilde{\bG}_{n})^{d}}\leq\frac{1}{k}$.
The pull-back $v_{k}\in\cV_{\pot}^{2}(\Omega)$ of $\tilde{v}_{k}$
satisfies $\bign\overline{w}_{i}-v_{k}\bign_{L^{2}(\bG_{n}^{,})^{d}}\leq\frac{1}{k}$
and hence $\overline{w}_{i}\in\cV_{\pot}^{2}(\bG_{n}\vert\S(\Rd))$.
We observe $\big(\tilde{w}_{i}^{(n)}+e_{i}\big)\I_{\tilde{\bG}_{n}}\in L_{\sol}^{2}(\S^{(n)})$
with the pull-back $\big(\overline{w}_{i}+e_{i}\big)\I_{\bG_{n}}\in L_{\sol}^{2}(\S(\Rd))$.
This implies $\overline{w}_{i}=w_{i}^{(n)}$. \Eqref{Pull-Back-Integral-Equation}
yields \Eqref{thm:pull-back-cell-solution}.
\end{proof}

\section{Proof of Lemma \ref{lem:Existence-of-Solution-eps} and \ref{lem:General-Homogenized-Limit}
\label{sec:Proof-of-Lemma}}

We first collect all the tools needed to prove the homogenization
result for minimally smooth domains (\lemref{General-Homogenized-Limit}).

\subsection{Extensions and traces for thinned point clouds}
\begin{thm}[Extending beyond holes and trace operator]
\label{thm:extension-trace-operators}\ \\
There exists a constant $C>0$ depending only on $n\in\N$ and $M_{0}>1$
such that the following holds:\\
Assume that $Q\subset\Rd$ is a bounded Lipschitz domain with Lipschitz
constant $M_{0}$, $\x\in\th\S(\Rd)$ and $\x_{Q}\subset\x$ such
that for every $x\in\x_{Q}$ it holds $\Ball{2r}x\subset Q$. Then
there exists an extension operator 
\[
\cU_{\x_{Q}}:\,W^{1,2}(Q\setminus\bm\x_{Q})\to W^{1,2}(Q)
\]
 such that $(\cU_{\x_{Q}}u)_{\vert Q\setminus\bm\x_{Q}}=u$ and 
\begin{equation}
\norm{\cU_{\x_{Q}}u}_{L^{2}(Q)}\leq C\norm u_{L^{2}(Q\setminus\bm\x_{Q})}\,,\qquad\norm{\nabla\cU_{\x_{Q}}u}_{L^{2}(Q)}\leq C\norm{\nabla u}_{L^{2}(Q\setminus\bm\x_{Q})}\,.\label{eq:thm:extension-trace-operators-1}
\end{equation}
Furthermore, there exists a trace operator 
\[
\cT_{\x_{Q}}:\,W^{1,2}(Q\setminus\bm\x_{Q})\to L^{2}(\partial\bm\x_{Q}):=L^{2}(\partial\bm\x_{Q},\,\H^{d-1})
\]
such that $\cT_{\x_{Q}}u=u|_{\partial\bm\x_{Q}}$ for every $u\in C_{c}^{1}(Q)$
and 
\begin{equation}
\bign\cT_{\x_{Q}}u\bign_{L^{2}(\partial\bm\x_{Q})}\leq C\left(\norm u_{L^{2}(Q\setminus\bm(\x_{Q}))}+\norm{\nabla u}_{L^{2}(Q\setminus\bm\x_{Q})}\right)\,.\label{eq:thm:extension-trace-operators-2}
\end{equation}
\end{thm}

\begin{proof}
For every $\x_{Q}\subset\x$ with $\x\in\th\S(\Rd)$, the set $Q\setminus\bm\x_{Q}$
is minimally smooth with $\delta=\min\left\{ \frac{1}{n},r\right\} $
and , $M=\max\left\{ \sqrt{2nr},M_{0}\right\} $. Furthermore, the
connected components of $\bm\x_{Q}$ have diameter less than $2nr$.
The existence of $\cU_{\x_{Q}}$ satisfying (\ref{eq:thm:extension-trace-operators-1})
follows from \cite[Proposition 3.3]{guillen2015quasistatic}. The
existence of $\cT_{\x_{Q}}$ satisfying (\ref{eq:thm:extension-trace-operators-2})
is provided in \cite{heida2021stochastic}.
\end{proof}

\subsection{Stochastic two-scale convergence}
\begin{defn}[Stationary and ergodic random measures]
\label{def:random-measure2}\ \\
A random measure $\mu_{\bullet}:\ \Omega\to\M(\R^{d})$ with underlying
dynamical system $(\Omega,\F,\p,\tau)$ is called \emph{stationary}
iff
\[
\mu_{\tau_{x}\omega}(A)=\mu_{\omega}(A+x)
\]
for every measurable $A\subset\R^{d}$, $x\in\R^{d}$ and $\p$-almost
every $\omega\in\Omega$. $\mu_{\bullet}$ is called \emph{ergodic}
iff it is stationary and $\tau$ is ergodic. 

\defref{random-measure2} is compatible with \defref{Random-Measure}
given in Section \ref{sec:Setting-Definitions-Main-Result} by considering
the canonical underlying probability space $(\M(\R^{d}),\B(\M(\R^{d}),\P_{\mu},\tau)$
with $\P_{\mu}$ being the distribution of $\mu$.
\end{defn}

\begin{thm}[Palm theorem (for finite intensity) \cite{Mecke1967}]
\label{thm:-Palm-Theorem} \ \\
Let $\mu_{\bullet}$ be a stationary random measure with underlying
dynamical system $(\Omega,\F,\p,\tau)$ of finite intensity $\lambda(\mu_{\bullet})$.\\
Then, there exists a unique finite measure $\mu_{\p}$ on $(\Omega,\F)$
such that for every $g:\ \R^{d}\times\Omega\to\R$ measurable and
either $g\geq0$ or $g\in L^{1}\left(\R^{d}\times\Omega,\,\L^{d}\otimes\mu_{\p}\right)$:
\[
\int_{\Omega}\int_{\R^{d}}g(x,\tau_{x}\omega)\dx\mu_{\omega}(x)\dx\p(\omega)=\int_{\R^{d}}\int_{\Omega}g(x,\omega)\dx\mu_{\p}(\omega)\dx x.
\]
For arbitrary $f\in L^{1}(\R^{d})$ with $\int_{\R^{d}}f\dx x=1$,
we have that
\[
\mu_{\p}(A)=\int_{\Omega}\int_{\R^{d}}f(x)\I_{A}(\tau_{x}\omega)\dx\mu_{\omega}(x)\dx\p(\omega)\,,
\]
in particular $\mu_{\p}(\Omega)=\lambda(\mu)$. Furthermore, for every
$\phi\in C_{c}(\Rd)$ and $g\in L^{1}(\Omega;\mupalm)$ the ergodic
limit 
\begin{equation}
\lim_{\eps\to0}\int_{\Rd}\phi(x)\,g\of{\tau_{\frac{x}{\eps}}\omega}\,\d x=\int_{\Rd}\int_{\Omega}\phi(x)\,g(\tilde{\omega})\,\d\mupalm(\tilde{\omega})\,\d x\label{eq:thm:-Palm-Theorem-ergodic}
\end{equation}
holds for $\p$-almost every $\omega$. We call $\mu_{\p}$ the Palm
measure of $\mu_{\bullet}$.
\end{thm}

For the rest of this subsection, we use the following assumptions.
\begin{assumption}
\label{assu:TS-assumption}\ \\
$\Omega$ is a compact metric space with a probability measure $\p$
and continuous dynamical system $\left(\tau_{x}\right)_{x\in\Rd}$.
Furthermore, $\mu_{\bullet}:\ \Omega\to\M(\R^{d})$ is a stationary
ergodic random measure with Palm measure $\mupalm$. We define $\muomegaeps(A):=\eps^{d}\muomega(\eps^{-1}A)$.
\end{assumption}

According to \cite{Zhikov2006} (by an application of \Eqref{thm:-Palm-Theorem-ergodic})
almost every $\omega\in\Omega$ is \emph{typical}, i.e. for such an
$\omega$, it holds for every $\phi\in C(\Omega)$ that
\[
\lim_{\eps\to0}\int_{Q}\phi\of{\tau_{\frac{x}{\eps}}\omega}\,\d x=\int_{\Omega}\phi\,\d\p\,.
\]

\begin{defn}[Two-scale convergence]
\label{def:ts-compact}\ \\
Let Assumption \ref{assu:TS-assumption} hold and let $\omega\in\Omega$
be typical. Let $\left(\ue\right)_{\eps>0}$ be a sequence $\ue\in L^{2}(Q,\muomegaeps)$
and let $u\in L^{2}(Q;L^{2}(\Omega,\mupalm))$ such that 
\[
\sup_{\eps>0}\norm{\ue}_{L^{2}(Q,\muomegaeps)}<\infty\,,
\]
and such that for every $\varphi\in C_{c}^{\infty}(Q)$, $\psi\in C(\Omega)$
\begin{equation}
\lim_{\eps\to0}\int_{Q}\ue(x)\varphi(x)\psi\of{\tau_{\frac{x}{\eps}}\omega}\,\d\muomegaeps(x)=\int_{Q}\int_{\Omega}u(x,\tilde{\omega})\varphi(x)\psi\of{\tilde{\omega}}\,\d\mupalm(\tilde{\omega})\,\d x\,.\label{eq:def:ts-compact}
\end{equation}
Then $\ue$ is said to be (weakly) \emph{two-scale convergent} to
$u$, written $\ue\tsweakto u$.
\end{defn}

\begin{rem}[Extending the space of test functions]
\label{rem:extended-test-functions}\ 
\begin{itemize}
\item If $\chi\in L^{\infty}(\Omega,\mupalm)$, then we can extend the class
of test functions from $\psi\in C(\Omega)$ to $\chi\psi$ since $\chi(\tau_{\frac{\cdot}{\eps}}\omega)\,\d\muomegaeps$
is again a random measure with Palm measure $\chi\d\mupalm$. 
\item Using a standard approximation argument, we can extend the class of
test functions from $\varphi\in C_{c}^{\infty}(Q)$ to $\varphi\in L^{2}(Q)$,
provided $\mu_{\omega}$ is uniformly continuous w.r.t the Lebesgue
measure. Then, strong $L^{2}(Q)$-convergence implies two-scale convergence
for $\mu_{\omega}\equiv\L^{d}$.
\end{itemize}
\end{rem}

\begin{lem}[{\cite[Lemma 5.1]{Zhikov2006}}]
\label{lem:Existence-ts-lim}\ \\
 Let Assumption \ref{assu:TS-assumption} hold. Let $1<p\leq\infty$,
$\omega\in\Omega$ be typical and $\ue\in L^{p}(Q,\muomegaeps)$ be
a sequence of functions such that $\norm{\ue}_{L^{p}(Q,\muomegaeps)}\leq C$
for some $C>0$ independent of $\eps$. Then there exists a subsequence
of $(u^{\eps'})_{\eps'\to0}$ and $u\in L^{p}(Q;L^{p}(\Omega,\mupalm))$
such that $u^{\eps'}\stackrel{2s}{\weakto}u$ and 
\begin{equation}
\norm u_{L^{p}(Q;L^{p}(\Omega,\mupalm))}\leq\liminf_{\eps'\to0}\norm{u^{\eps'}}_{L^{p}(Q,\muomegaeps)}\,.\label{eq:two-scale-limit-estimate}
\end{equation}
\end{lem}

\begin{thm}[Two-scale convergence in $W^{1,2}(Q)$ \cite{Zhikov2006}]
\label{thm:sto-conver-grad}\ \\
Under Assumption \ref{assu:TS-assumption}, for every typical $\omega\in\Omega$
the following holds: \\
If $u^{\eps}\in W^{1,2}(Q;\Rd)$ for all $\eps$ and if 
\[
\sup_{\eps>0}\big(\norm{u^{\eps}}_{L^{2}(Q)}+\norm{\nabla u^{\eps}}_{L^{2}(Q)}\big)<\infty\,,
\]
then there exists a $u\in W^{1,2}(Q)$ with $u^{\eps}\weakto u$ weakly
in $W^{1,2}(Q)$ and there exists $v\in L^{2}(Q;\cV_{\pot}^{2}(\Omega))$
such that $\nabla\ue\tsweakto\nabla u+v$ weakly in two scales.
\end{thm}

\subsection{\label{subsec:Two-scale-P}Two-scale convergence on perforated domains}

Due to Theorem \ref{thm:Fn-compact-and-nice}, the set $\th\S(\Rd)\subset\M(\Rd)$
is compact, hence the above two-scale convergence method can be applied
for the stationary ergodic point process $\X^{(n)}$ taking values
in $\th\S(\R^{d})$ only. To be more precise, we consider the compact
metric probability space $\Omega=\th\S(\Rd)$ and a random variable
$\X^{n}:\ \Omega\to\S(\R^{d})$ such that $\X^{n}$ and $\X^{(n)}$
have the same distribution. By the considerations made in Subsection
\ref{subsec:Pull-back}, they will both result in the same partial
differential equation.
\begin{thm}[Extension and trace estimates on $Q_{\x}^{\eps}$ for $\x\in\th\S(\R^{d})$]
\label{thm:Extension-Trace-Estimates-eps}\ \\
Let $Q\subset\Rd$ be a bounded domain, $n\in\N$ be fixed. Let $\X$
be an admissible point process with values in $\th\S(\R^{d})$. For
almost every realization $\x$ of $\X$, we have: \\
Let $Q_{\x}^{\eps}$ and $G_{\x}^{\eps}$ be defined according to
\defref{perforation-Q}.
\begin{enumerate}
\item There exists a $C>0$ depending only on $Q$ and $n$ and a family
of extension and trace operators 
\[
\cU_{\eps,\x}:\,W^{1,2}(Q_{\x}^{\eps})\to W^{1,2}(Q)\,,\qquad\cT_{\eps,\x}:\,W^{1,2}(Q_{\x}^{\eps})\to L^{2}(\partial G_{\x}^{\eps})
\]
such that for every $u\in W^{1,2}(Q_{\x}^{\eps})$ it holds
\begin{align*}
\norm{\cU_{\eps,\x}u}_{W^{1,2}(Q)} & \leq C\norm u_{W^{1,2}(Q_{\x}^{\eps})}\,,\\
\eps\norm{\cT_{\eps,\x}u}_{L^{2}(\partial G_{\x}^{\eps})}^{2} & \leq C\left(\norm u_{L^{2}(Q_{\x}^{\eps})}^{2}+\eps^{2}\norm{\nabla u}_{L^{2}(Q_{\x}^{\eps})}^{2}\right)\,.
\end{align*}
\item If $u^{\eps}\in W^{1,2}(Q_{\x}^{\eps})$ is a sequence satisfying
$\sup_{\eps}\norm{\ue}_{W^{1,2}(Q_{\x}^{\eps})}<\infty$, then there
exists a $u\in W^{1,2}(Q)$ and a subsequence still indexed by $\eps$
such that $\cU_{\eps,\x}\ue\weakto u$ weakly in $W^{1,2}(Q)$ and
there exists $v\in L^{2}(Q;L_{\pot}^{2}(\Omega))$ such that
\[
\nabla\cU_{\eps,\x}\ue\tsweakto\nabla u+v\,,\qquad\nabla\ue\tsweakto\I_{\bG_{n}}\left(\nabla u+v\right)\,,
\]
with the event $\mathbf{G}_{n}:=\left\{ \x\in\th\S(\R^{d})\,\vert\,o\notin\bm\x\right\} $.
Furthermore, for some $C>0$ depending only on $Q$ and $n$
\begin{equation}
\eps\norm{\cT_{\eps,\x}(\ue-u)}_{L^{2}(\partial G_{\x}^{\eps})}^{2}\leq C\left(\norm{\cU_{\eps,\x}\ue-u}_{L^{2}(Q)}^{2}+\eps^{2}\norm{\nabla\cU_{\eps,\x}\ue-\nabla u}_{L^{2}(Q)}^{2}\right)\,.\label{eq:Trace-to-0}
\end{equation}
\end{enumerate}
\end{thm}

\begin{proof}
1. follows from using \thmref{extension-trace-operators} on $\eps^{-1}G_{\x}^{\eps}$
and rescaling the inequalities (\ref{eq:thm:extension-trace-operators-1})--(\ref{eq:thm:extension-trace-operators-2}). 

2. is a bit more lengthy. The existence of a subsequence and $u\in W^{1,2}(Q)$
and $v\in L^{2}(Q;L_{\pot}^{2}(\Omega))$ such that $\cU_{\eps,\x}\ue\weakto u$
and $\nabla\cU_{\eps,\x}\ue\tsweakto\nabla u+v$ follows from \thmref{sto-conver-grad}.
We observe that $\I_{\mathbf{G}_{n}}(\tau_{x}\x)=\I_{\bm\x^{\complement}}(x)$
and $\I_{\bm\x^{\complement}}(\frac{x}{\eps})=\I_{\eps\bm\x^{\complement}}(x)$.
Therefore, $\I_{\eps\bm\x^{\complement}}\nabla\cU_{\eps,\x}\ue\tsweakto\I_{\mathbf{G}_{n}}\of{\nabla u+v}$
(\remref{extended-test-functions} (1)). Furthermore, we observe with
$Q_{n,r}^{\eps}:=\left\{ x\in Q:\;\dist(x,\partial Q)\leq\eps nr\right\} $
that $\left(\eps\bm\x\cap Q\right)\backslash G_{\x}^{\eps}\subset Q_{n,r}^{\eps}$
and
\[
\biga\I_{\eps\bm\x}-\I_{G_{\x}^{\eps}}\biga\leq\I_{Q_{n,r}^{\eps}}\xrightarrow{\eps\to0}0\text{ pointwise a.e.}
\]
Therefore, $\I_{\eps\bm\x}-\I_{G_{\x}^{\eps}}\to0$ strongly in $L^{p}(Q)$,
$p\in[1,\infty)$ and hence taking $\phi\in C(\Omega)$, $\psi\in C(\overline{Q})$,
we find 
\[
\int_{Q}\of{\I_{\eps\bm\x}-\I_{G_{\x}^{\eps}}}\nabla\cU_{\eps,\x}\ue\phi(\tau_{\frac{\bullet}{\eps}}\x)\psi\leq\norm{\I_{\eps\bm\x}-\I_{G_{\x}^{\eps}}}_{L^{2}(Q)}\norm{\nabla\cU_{\eps,\x}\ue}_{L^{2}(Q)}\norm{\phi}_{\infty}\norm{\psi}_{\infty}\to0\,.
\]
In particular, $\I_{\eps\bm\x}\nabla\cU_{\eps,\x}\ue$ and $\nabla u^{\eps}=\I_{G_{\x}^{\eps}}\nabla\cU_{\eps,\x}\ue$
have the same two-scale limit (\remref{extended-test-functions} (2))
\[
\nabla u^{\eps}\xrightharpoonup[\eps\to0]{2s}\I_{\mathbf{G}}\of{\nabla u+v}\,.
\]
Due to the absolutely bounded diameter of the connected components
of $\bm\x$, there exists a domain $B\supset Q$ big enough such that,
with the notation of Definition \ref{def:perforation-Q},
\[
Q\cap\eps\bm(J_{\eps}(\x,B))=Q\cap\eps\bm\x\quad\forall\eps\in(0,1)\,.
\]
Now let $\cU_{Q}:\,W^{1,2}(Q)\to W^{1,2}(B)$ be the canonical extension
operator satisfying 
\[
\norm{\cU_{Q}u}_{L^{2}(B)}\leq C\norm{\cU_{Q}u}_{L^{2}(Q)}\qquad\text{and}\qquad\norm{\nabla\cU_{Q}u}_{L^{2}(B)}\leq C\norm{\nabla u}_{L^{2}(Q)}\,.
\]
 Reapplying \thmref{extension-trace-operators} to the trace on $\eps^{-1}\left(B\backslash\eps\bm(J_{\eps}(\x,B))\right)$,
we find for some constant $C$ independent from $\eps$ and $\X$
but depending on $Q$, $B$ and $n$ and varying from line to line:
\begin{align*}
\eps\norm{\cT_{\eps,\x}(\ue-u)}_{L^{2}(\partial G_{\x}^{\eps})}^{2} & \leq\eps\norm{\cT_{\eps,\x}(\cU_{Q}\cU_{\eps,\x}\ue-\cU_{Q}u)}_{L^{2}(\eps\partial\bm(J_{\eps}(\x,B)))}^{2}\\
 & \leq C\left(\norm{\cU_{Q}\cU_{\eps,\x}\ue-\cU_{Q}u}_{L^{2}(B)}^{2}+\eps^{2}\norm{\nabla\cU_{Q}\cU_{\eps,\x}\ue-\nabla\cU_{Q}u}_{L^{2}(B)}^{2}\right)\\
 & \leq C\left(\norm{\cU_{\eps,\x}\ue-u}_{L^{2}(Q)}^{2}+\eps^{2}\norm{\nabla\cU_{\eps,\x}\ue-\nabla u}_{L^{2}(Q)}^{2}\right)\,.
\end{align*}
\end{proof}

\subsection{Existence of solution on perforated domains (Lemma \ref{lem:Existence-of-Solution-eps})}

Due to the perforations, $\partial_{t}\ue$ cannot be embedded in
a common space in a convenient way for the application of the \foreignlanguage{british}{Aubin--Lions}
theorem. Hence we use the following general characterization of compact
sets. 
\begin{thm}[{Characterization of compact sets in $L^{p}(I;V)$ \cite[Theorem 1]{simon1986compact}}]
\label{thm:Simon-Compactness}\ \\
Let $V$ be a Banach space, $p\in[1,\infty)$ and $\Lambda\subset L^{p}(I;V)$.
$\Lambda$ is relatively compact in $L^{p}(I;V)$ if and only if
\begin{align}
 & \left\{ \int_{t_{1}}^{t_{2}}v(t)\dx t\,\vert\,v\in\Lambda\right\} \text{ is relatively compact in }V\ \forall\,0<t_{1}<t_{2}<T,\label{eq:Simon-1}\\
 & \sup_{v\in\Phi}\bign\mathfrak{s}_{h}[v]-v\bign_{L^{p}(0,T-h;V)}\to0\text{ as }h\to0\,,\label{eq:Simon-2}
\end{align}
where $\mathfrak{s}_{h}[v(\,\cdot\,)]:=v(\,\cdot\,+h)$ is the shift
by $h\in\left(0,\,T\right)$.
\end{thm}

We can now establish the existence of a solution for fixed $\eps>0$
to our partial differential equation. 
\begin{thm}[Existence of solution on perforated domains and a priori estimate]
\label{thm:Solution-Perforated-Domain}~\\
Let $\x\in\th\S(\R^{d})$. Under \assuref{PDE-Parameters} and with
$Q_{\x}^{\eps}$ as defined in \defref{perforation-Q}, we have:\\
There exists a solution $\uexn\in L^{2}(I;W^{1,2}(Q_{\x}^{\eps}))$
with generalized time derivative $\del_{t}\uexn\in L^{2}(I;W^{1,2}(Q_{\x}^{\eps})^{*})$
to \Eqref{System-in-Eps}, i.e.
\begin{equation}
\begin{aligned}\del_{t}\uexn-\nabla\cdot\left(A(\uexn)\,\nabla\uexn\right) & =f &  & \text{in }I\times Q_{\x}^{\eps}\\
A(\uexn)\,\nabla\uexn\cdot\nu & =0 &  & \text{on }I\times\del Q\\
A(\uexn)\,\nabla\uexn\cdot\nu & =\eps h(\uexn) &  & \text{on }I\times\del Q_{\x}^{\eps}\backslash\del Q\\
\uexn(0,x) & =u_{0}(x) &  & \text{in }Q_{\x}^{\eps}\,,
\end{aligned}
\label{eq:System-Perforated-u-eps}
\end{equation}
which satisfies the following a priori estimates for $\eps$ small
enough
\begin{align}
\esssup_{t\in I}\bign\uexn(t)\bign_{L^{2}(Q_{\x}^{\eps})}^{2} & \leq\exp(C_{1})\big[\bign u_{0}\bign_{L^{2}(Q)}^{2}+C_{2}\big]\nonumber \\
\bign\nabla\uexn\bign_{L^{2}(I;L^{2}(Q_{\x}^{\eps}))}^{2} & \leq\frac{1}{\inf(A)}\big(1+C_{1}\exp(C_{1})\big)\big[\bign u_{0}\bign_{L^{2}(Q)}^{2}+C_{2}\big]\label{eq:A-Priori-Eps}\\
\bign\del_{t}\uexn\bign_{L^{2}(I;W^{1,2}(Q_{\x}^{\eps})^{*})}^{2} & \leq\tilde{C}\,,\nonumber 
\end{align}
where
\[
C_{1}:=T(1+3CL_{h})\qquad\text{and}\qquad C_{2}:=TL_{h}h(0)^{2}\L^{d}(Q)+\bign f\bign_{L^{1}(I;L^{2}(Q))}^{2}\,,
\]
$C$ is from \thmref{Extension-Trace-Estimates-eps} depending only
on $Q$ and $n$ and where $\tilde{C}>0$ is independent of $\eps$.
\end{thm}

\begin{proof}
We will only sketch the proof. There are 3 main steps: Deriving a
priori estimates, existence of Galerkin solutions and the limit passing.

1. Testing \Eqref{System-Perforated-u-eps} with $\uexn$ and using
\[
\langle\del_{t}\uexn,\,\uexn\rangle_{W^{1,2}(Q_{\x}^{\eps})^{*},W^{1,2}(Q_{\x}^{\eps})}=\frac{1}{2}\frac{\dx}{\dx t}\bign\uexn\bign_{L^{2}(Q_{\x}^{\eps})}^{2}
\]
 yields
\begin{align*}
\frac{1}{2}\frac{\dx}{\dx t}\bign\uexn\bign_{L^{2}(Q_{\x}^{\eps})}^{2}+A(\uexn)\,\bign\nabla\uexn\bign_{L^{2}(Q_{\x}^{\eps})}^{2}-\eps\big(h(\uexn),\,\uexn\big)_{L^{2}(\del G_{\x}^{\eps})} & =\big(f,\,\uexn\big){}_{L^{2}(Q_{\x}^{\eps})}\,.
\end{align*}
The a priori estimate then follows from the Gronwall inequality and
the trace estimate in \thmref{Extension-Trace-Estimates-eps}.\\
For the a priori estimate in $\del_{t}\uexn$, one simply uses\textcolor{teal}{}
\[
\langle\del_{t}\uexn,\,\varphi\rangle=\big(A(\uexn)\,\nabla\uexn,\,\nabla\varphi\big){}_{L^{2}(Q_{\x}^{\eps})}+\eps\big(h(\uexn),\,\varphi\big){}_{L^{2}(\del G_{\x}^{\eps})}+\big(f,\,\varphi\big){}_{L^{2}(Q_{\x}^{\eps})}\,.
\]

2. Let $(V_{m})_{m\in\N}$ be a family of finite-dimensional vector
spaces, $V_{m}\nearrow W^{1,2}(Q_{\x}^{\eps})$. One can show that
solutions to \Eqref{System-in-Eps} exist in $V_{m}$, e.g. via fixed
point arguments. These solutions $u_{(m)}^{\eps}$ also satisfy the
a priori estimate in \Eqref{A-Priori-Eps} and 
\begin{equation}
\sup_{m\in\N}\bign\del_{t}u_{(m)}^{\eps}\bign_{L^{2}(I;V_{m}^{*})}<\infty\,.\label{eq:Galerkin-A-Priori-Time}
\end{equation}

3. The a priori estimates yield a $L^{2}(I;L^{2}(Q_{\x}^{\eps}))$-weakly
convergent subsequence to some $\uexn\in L^{2}(I;L^{2}(Q_{\x}^{\eps}))$.
\thmref{Simon-Compactness} and \Eqref{Galerkin-A-Priori-Time} imply
pre-compactness of $(u_{(m)}^{\eps})_{m\in\N}\subset L^{2}(I;L^{2}(Q_{\x}^{\eps}))$
as well as pre-compactness of $(\cT_{\eps,\x}u_{(m)}^{\eps})_{m\in\N}\subset L^{2}(I;L^{2}(\del G_{\x}^{\eps}))$,
see \remref{Procedure-Simon}. Testing with functions in $L^{2}(I;V_{m})$
and passing to the limit $m\to\infty$ finishes the proof since $\bigcup_{m\in\N}V_{m}$
is dense in $W^{1,2}(Q_{\x}^{\eps})$.
\end{proof}
\begin{rem}[Procedure of Simon's theorem]
\label{rem:Procedure-Simon}\ \\
We will use Simon's theorem (\thmref{Simon-Compactness}) on multiple
occasions. The general procedure will always be the same. We will
exemplary prove the following result:\\
Let $I=[0,T]$, $U\subset\R^{d}$ be some bounded Lipschitz-domain
and $\cT:\ W^{1,2}(U)\to L^{2}(\del U)$ the trace operator. For each
$k\in\N$, let $u_{k}\in L^{2}(I;W^{1,2}(U))$ with generalized time-derivative
$\del_{t}u_{k}\in L^{2}(I;W^{1,2}(U)^{*})$ via $W^{1,2}(U)\hookrightarrow L^{2}(U)\hookrightarrow W^{1,2}(U)^{*}$.
Assume that
\[
C:=\sup_{k\in\N}\bign u_{k}\bign_{L^{2}(I;W^{1,2}(U))}<\infty\qquad\text{and}\qquad\tilde{C}:=\sup_{k\in\N}\bign\del_{t}u_{k}\bign_{L^{2}(I;V_{k}^{*})}<\infty
\]
for either the situation that $W^{1,2}(U)\subset V_{k}\subset L^{2}(U)$
with $\norm{\cdot}_{V_{k}}\leq\norm{\cdot}_{W^{1,2}(U)}$ and uniformly
continuous injective maps $\cU_{k}:\,V_{k}\to W^{1,2}(U)$ or for
the situation that $V_{k}\subset W^{1,2}(U)$. We further claim $u_{k}(t)\in V_{k}$
for almost every $t\in I$. Then, 
\[
(\cU_{k}u_{k})_{k\in\N}\subset L^{2}(U)\text{ resp. }(u_{k})_{k\in\N}\subset L^{2}(U)\quad\text{and }\quad(\cT u_{k})_{k\in\N}\subset L^{2}(\del U)
\]
are relatively compact.
\end{rem}

\begin{proof}[Exemplary proof for the procedure of Simon's theorem]
 We need to show Condition (\ref{eq:Simon-1}) and Condition (\ref{eq:Simon-2})
from \thmref{Simon-Compactness}.
\begin{enumerate}
\item Condition (\ref{eq:Simon-1}) usually relies on compactness results
for the stationary setting. Since
\[
\sup_{k\in\N}\bign\int_{t_{1}}^{t_{2}}u_{k}\dx t\bign_{W^{1,2}(Q)}\leq\sup_{k\in\N}\sqrt{T}\bign u_{k}\bign_{L^{2}(I;W^{1,2}(U))}<\infty,
\]
compactness of $\cT$ yields pre-compactness of $\big(\int_{t_{1}}^{t_{2}}\cT u_{k}\dx t\big)_{k\in\N}=\big(\cT\int_{t_{1}}^{t_{2}}u_{k}\dx t\big)_{k\in\N}\subset L^{2}(\del U)$,
so we have shown Condition (\ref{eq:Simon-1}).
\item Condition (\ref{eq:Simon-2}) will additionally require some a-priori-estimate
on $\del_{t}u_{k}$. We have
\[
u_{k}(t_{2})=u_{k}(t_{1})+\int_{t_{1}}^{t_{2}}\del_{t}u_{k}\dx s
\]
as elements of $W^{1,2}(U)^{*}$. Using the Cauchy--Schwarz inequality
twice, we get for $h\in(0,T)$:
\begin{align*}
\bign\mathfrak{s}_{h}[u_{k}]-u_{k} & \bign_{L^{2}((0,T-h);L^{2}(U))}^{2}=\int_{0}^{T-h}\big(u_{k}(t+h)-u_{k}(t),\,u_{k}(t+h)-u_{k}(t)\big)_{L^{2}(U)}\dx t\\
 & =\int_{0}^{T-h}\bigl\int_{t}^{t+h}\del_{t}u_{k}(s)\dx s,\,u_{k}(t+h)-u_{k}(t)\bigr_{W^{1,2}(U)^{*},W^{1,2}(U)}\dx t\\
 & \leq\int_{0}^{T-h}\bign\int_{t}^{t+h}\del_{t}u_{k}(s)\dx s\bign_{L^{2}(V_{k}^{*})}\,\bign\cU_{k}u_{k}(t+h)-\cU_{k}u_{k}(t)\bign_{W^{1,2}(U)}\dx t\\
 & \leq h\bign\del_{t}u_{k}\bign_{L^{2}(I;V_{k}^{*})}\,2\bign u_{k}\bign_{L^{2}(I;W^{1,2}(U))}\leq2hC\tilde{C}\,.
\end{align*}
Compactness of $\cT$ implies that for every $\delta>0$, there exists
a $C_{\delta}>0$ such that
\[
\bign\cT v\bign_{L^{2}(\del U)}^{2}\leq C_{\delta}\bign v\bign_{L^{2}(U)}^{2}+\delta\bign\nabla v\bign_{L^{2}(U)}^{2}\ \forall v\in W^{1,2}(U)\,.
\]
Therefore,
\begin{align*}
\bign\mathfrak{s}_{h}[\cT u_{k}] & -\cT u_{k}\bign_{L^{2}((0,T-h);L^{2}(\del U))}^{2}=\bign\cT[\mathfrak{s}_{h}u_{k}-u_{k}]\bign_{L^{2}((0,T-h);L^{2}(\del U))}^{2}\\
\leq & C_{\delta}\bign\mathfrak{s}_{h}u_{k}-u_{k}\bign_{L^{2}(\del U)}^{2}+\delta\bign\nabla\mathfrak{s}_{h}u_{k}-\nabla u_{k}\bign_{L^{2}(U)}^{2}\\
\leq & 2hC_{\delta}C\tilde{C}+2\delta\tilde{C}\,.
\end{align*}
The estimate is independent of the chosen $u_{k}$, Condition (\ref{eq:Simon-2})
holds.
\end{enumerate}
We have shown both conditions and conclude.
\end{proof}

\subsection{Homogenization for minimally smooth domains (Lemma \ref{lem:General-Homogenized-Limit})}

We can now pass to the limit $\eps\to0$ for the homogenized system.
Some extra care has to be taken since $Q_{\x}^{\eps}\neq Q\backslash\eps\bm\x$,
especially in the boundary term. However, we show that the difference
becomes negligible for the two-scale convergence as $\eps\to0$.
\begin{thm}[Homogenized system for $\bm\X^{(n)}$]
\label{thm:Homogenized-System-for-n}\ \\
Let $\X$ be a stationary ergodic point process with values in $\th\S(\R^{d})$.
Recall the surface measure $\mu_{\x}$ from \defref{Surface-Measure-of-BM}
\[
\mu_{\x}:=\H_{\llcorner\del\bm\x}^{d-1}\,.
\]
Under \assuref{PDE-Parameters}, we have for almost every realization
$\x$ of $\X$ and with $Q_{\x}^{\eps}$ as defined in \defref{perforation-Q}:
\\
Let $u^{\eps}$ be a solution to \Eqref{System-Perforated-u-eps}
and let $\cU_{\eps,\x}$ be given as in \thmref{Extension-Trace-Estimates-eps}.
There exists a $u_{n}\in L^{2}(I;W^{1,2}(Q))$ with generalized time
derivative $\del_{t}u_{n}\in L^{2}(I;W^{1,2}(Q)^{*})$ such that for
a subsequence
\begin{align*}
\cU_{\eps,\x}u^{\eps} & \xrightarrow[\eps\to0]{L^{2}(I;L^{2}(Q))}u_{n}\\
\del_{t}u^{\eps} & \xrightharpoonup[\eps\to0]{L^{2}(I;W^{1,2}(Q)^{*})}\P(\bG_{n})\del_{t}u_{n}
\end{align*}
and $u_{n}$ is a (not necessarily unique) solution to 
\begin{align}
\P(\bG_{n})\del_{t}u_{n}-\nabla\cdot\left(A(u_{n})\A\,\nabla u_{n}\right)-\lambda(\mu_{\X})h(u_{n}) & =\P(\bG_{n})f &  & \text{in }I\times Q\nonumber \\
A(u_{n})\A^{(n)}\,\nabla u_{n}\cdot\nu & =0 &  & \text{on }I\times\del Q\label{eq:System-u_n}\\
u_{n}(0,x) & =\P(\bG_{n})u_{0}(x) &  & \text{in }Q\nonumber 
\end{align}
with $\A^{(n)}$ being the effective conductivity based on the event
$\bG_{n}=\big\{\x\in\th\S(\R^{d})\,\vert\,o\notin\bm\x\big\}$ defined
in \defref{Effective-Conductivity}. Furthermore, $u_{n}$ satisfies
the following a priori estimates
\begin{align*}
\esssup_{t\in I}\bign u_{n}(t)\bign_{L^{2}(Q)}^{2} & \leq\exp(C_{1}^{(n)})\big[\bign u_{0}\bign_{L^{2}(Q)}^{2}+C_{2}^{(n)}\big]\\
\bign\nabla u_{n}\bign_{L^{2}(I;L^{2}(Q))}^{2} & \leq\frac{\P(\bG_{n})}{2\alpha_{\A^{(n)}}\inf(A)}\big(1+C_{1}^{(n)}\exp(C_{1}^{(n)})\big)\big[\bign u_{0}\bign_{L^{2}(Q)}^{2}+C_{2}^{(n)}\big]
\end{align*}
for
\begin{align*}
C_{1}^{(n)} & :=T\big(1+\frac{\lambda(\mu_{\X})}{\P(\bG_{n})}(1+2L_{h})\big)\\
C_{2}^{(n)} & :=\bign f\bign_{L^{2}(I;L^{2}(Q))}^{2}+2T\frac{\lambda(\mu_{\X})}{\P(\bG_{n})}\biga h(0)\biga^{2}.
\end{align*}
\end{thm}

\begin{proof}
The a priori estimates in \Eqref{A-Priori-Eps} and \thmref{Extension-Trace-Estimates-eps}
tell us that
\begin{align*}
\cU_{\eps,\x}u^{\eps} & \tsweakto u_{n}\,, & \nabla\cU_{\eps,\x}u^{\eps} & \tsweakto\nabla u_{n}+v\,,\\
u^{\eps} & \tsweakto\I_{\bG_{n}}u_{n}\,, & \nabla u^{\eps} & \tsweakto\I_{\bG_{n}}\nabla u_{n}+v\,,
\end{align*}
for some $u_{n}\in L^{2}(I;W^{1,2}(Q))$ and $v\in L^{2}(I;L^{2}(Q;L_{\pot}^{2}(\Omega)))$
where the two-scale convergence is with respect to the Lebesgue measure
$\L^{d}$. The uniform bound for $\del_{t}u^{\eps}$ in \Eqref{A-Priori-Eps}
together with \thmref{Simon-Compactness} yields (for yet another
subsequence) 
\begin{equation}
\cU_{\eps,\x}\ue\xrightarrow[\eps\to0]{L^{2}(I;L^{2}(Q))}u_{n}\,,\label{eq:Strong-L2-Convergence-Eps}
\end{equation}
compare to, e.g., \remref{Procedure-Simon}.

For $\varphi_{1},\varphi_{2}\in C^{1}([0,T]\times\overline{Q})$ with
$\varphi_{1}(T,\cdot)=0$ and $\psi\in H^{1}(\Omega)$, we use $\varphi^{\eps}(t,x):=\varphi_{1}(t,x)+\eps\varphi_{2}(t,x)\psi(\tau_{\frac{x}{\eps}}\x)$
as a test function and pass to the limit using two-scale convergence.
Furthermore, we use 
\begin{equation}
A(u^{\eps})\nabla u^{\eps}\tsweakto A(u_{n})\I_{\bG_{n}}\left(\nabla u_{n}+v\right)\qquad\text{and}\qquad h(u^{\eps})\xrightharpoonup{2s,\,\mu_{\X}}h(u_{n})\,,\label{eq:thm:Homogenized-System-for-n-help-1}
\end{equation}
which we will prove below. We then obtain the two equations 
\begin{align*}
{-}\int_{0}^{T}\int_{Q}u_{n}\partial_{t}\varphi_{1}+\int_{Q}u_{0}\varphi_{1}+\int_{0}^{T}\int_{Q}\int_{\bG_{n}}\nabla\varphi_{1}\cdot A(u_{n})\left(\nabla u_{n}+v\right)\\
+\int_{0}^{T}\int_{Q}h(u_{n})\varphi_{1}\int_{\Omega}\d\mupalm & =\int_{0}^{T}\int_{Q}f\varphi_{1}\,,\\
\int_{0}^{T}\int_{Q}\int_{\bG_{n}}\varphi_{2}\nabla_{\omega}\psi\cdot A(u_{n})\left(\nabla u_{n}+v\right) & =0\,.
\end{align*}
The second equation holds true for every choice of $\varphi_{2}$
and $\psi$ as above if we make the standard ansatz $v=\sum_{i=1}^{d}\partial_{i}u_{n}w_{i}^{(n)}$,
where $w_{i}^{(n)}$ are the cell solutions from \defref{Cell-Solutions}
for $\Omega=\th\S(\Rd)$ and $\p=\P$ being the distribution of $\X$.
Plugging this information into the first equation yields \Eqref{System-u_n}.
The a priori estimate follows from testing \Eqref{System-u_n} with
$u_{n}$ and the Gronwall inequality (see e.g. the proof of \thmref{Solution-Perforated-Domain}).\textcolor{red}{{}
}It only remains to prove \Eqref{thm:Homogenized-System-for-n-help-1}.

Now, we show the first part of \Eqref{thm:Homogenized-System-for-n-help-1}.
By \remref{extended-test-functions}, we know that 
\[
A(u_{n})\nabla u^{\eps}\xrightharpoonup[\eps\to0]{2s}\I_{\bG_{n}}A(u_{n})\of{\nabla u_{n}+v}\,.
\]
Using dominated convergence and \Eqref{Strong-L2-Convergence-Eps}
yields a subsequence such that $A(\cU_{\eps,\x}u^{\eps})\to A(u_{n})$
in $L^{p}(0,T;L^{p}(Q))$ for every $1\leq p<\infty$. Using test
functions $\phi\in C(\overline{Q})$ and $\psi\in C(\Omega)$ we observe
that $\left(A(\cU_{\eps,\x}u^{\eps})-A(u_{n})\right)\nabla u^{\eps}\tsweakto0$,
so $A(u^{\eps})\tsweakto\I_{\bG_{n}}A(u_{n})\of{\nabla u_{n}+v}$. 

The second part of \Eqref{thm:Homogenized-System-for-n-help-1} is
more difficult. Given $\varphi\in C^{1}(\overline{Q})$ and $\psi\in C(\Omega)$,
we set $\psi^{\eps,\x}(x):=\psi(\tau_{\frac{x}{\eps}}\x)$ and find
\begin{align*}
 & \left|\eps\int_{\del G_{\x}^{\eps}}h\big(\ue(x)\big)\varphi(x)\psi^{\eps,\x}(x)\d\H^{d-1}(x)-\int_{Q}\int_{\Omega}h\big(u_{n}(x)\big)\,\varphi(x)\,\psi(\x)\,\d\mupalm(\x)\,\d x\right|\\
 & \leq\quad\left|\eps\int_{\del G_{\x}^{\eps}}h\big(\ue(x)\big)\varphi(x)\psi^{\eps,\x}(x)\d\H^{d-1}(x)-\eps\int_{\del G_{\x}^{\eps}}h\big(u_{n}(x)\big)\varphi(x)\psi^{\eps,\x}(x)\d\H^{d-1}(x)\right|\\
 & \quad+\left|\eps\int_{\del G_{\x}^{\eps}}h\big(u_{n}(x)\big)\varphi(x)\psi^{\eps,\x}(x)\d\H^{d-1}(x)-\eps\int_{Q\cap\eps\partial\bm\x}h\big(u_{n}(x)\big)\varphi(x)\psi^{\eps,\x}(x)\d\H^{d-1}(x)\right|\\
 & \quad+\left|\eps\int_{Q\cap\eps\partial\bm\x}h\big(u_{n}(x)\big)\varphi(x)\psi^{\eps,\x}(x)\d\H^{d-1}(x)-\int_{Q}\int_{\Omega}h\big(u_{n}(x)\big)\,\varphi(x)\,\psi(\x)\,\d\mupalm(\x)\,\d x\right|\,.
\end{align*}
We will show that all these terms go to $0$ as $\eps\to0$. Due to
the Lipschitz continuity of $h$ and Stampaccias lemma we find
\begin{align*}
\norm{\nabla h(u^{\eps})}_{L^{2}(Q_{\x}^{\eps})} & \leq\norm h_{C^{0,1}}\norm{\nabla u^{\eps}}_{L^{2}(Q_{\x}^{\eps})}\,, & \norm{h(u^{\eps})}_{L^{2}(Q_{\x}^{\eps})} & \leq\norm h_{C^{0,1}}\left(\norm{u^{\eps}}_{L^{2}(Q_{\x}^{\eps})}+1\right)\,,\\
\norm{\nabla h(u_{n})}_{L^{2}(Q)} & \leq\norm h_{C^{0,1}}\norm{\nabla u_{n}}_{L^{2}(Q)}\,, & \norm{h(u_{n})}_{L^{2}(Q)} & \leq\norm h_{C^{0,1}}\left(\norm{u_{n}}_{L^{2}(Q)}+1\right)\,.
\end{align*}
Furthermore, $\cU_{\eps,\x}\ue\to u_{n}$ strongly in $L^{2}(I;L^{2}(Q))$
and weakly in $L^{2}(I;W^{1,2}(Q))$ implies $h\left(\cU_{\eps,\x}\ue\right)\to h(u_{n})$
in the same topologies. $G_{\x}^{\eps}$ as in \defref{perforation-Q}
fulfills $\del G_{\x}^{\eps}=\del Q_{\x}^{\eps}\backslash\del Q$.
\Eqref{Trace-to-0} together with the strong convergence of $\cU_{\eps,\x}u^{\eps}$
tell us
\[
\eps\norm{h\big(\cT_{\eps,\x}\ue\big)-h\big(\cT_{\eps,\x}u_{n}\big)}_{L^{2}(I;L^{2}(\partial G_{\x}^{\eps}))}^{2}\leq L_{h}\eps\norm{\cT_{\eps,\x}(\ue-u_{n})}_{L^{2}(I;L^{2}(\partial G_{\x}^{\eps}))}^{2}\to0\,,
\]
which already shows convergence in the first summand. Similar considerations
to the proof of \Eqref{Trace-to-0} tell us that $\cT_{\eps,\x}:\ W^{1,2}(Q_{\x}^{\eps})\to L^{2}(\eps\del\bm\x\cap Q)$
is a bounded linear operator, so we can consider the trace not only
on $\partial G_{\x}^{\eps}$ but even for clusters close to the boundary.
We have, with $C>0$ changing from line to line but independent of
$\eps$,
\begin{align*}
\biga\eps\int_{I}\int_{Q} & \big(\I_{\eps\del\bm\x}-\I_{\del G_{\x}^{\eps}}\big)h(u^{\eps})\varphi\psi^{\eps,\x}\dx\H^{d-1}\dx t\biga^{2}\\
 & \leq C\eps\bign h(u^{\eps})\bign_{L^{2}(I;L^{2}((\eps\del\bm\x)\backslash\del G_{\x}^{\eps}))}^{2}\cdot\eps\bign1\bign_{L^{2}(I;L^{2}((\eps\del\bm\x)\backslash\del G_{\x}^{\eps}))}^{2}\\
 & \leq C\big\{\bign\I_{Q_{n,r}^{\eps}}h(u^{\eps})\bign_{L^{2}(I;L^{2}(Q))}^{2}+\eps\bign\I_{Q_{n,r}^{\eps}}\nabla h(u^{\eps})\bign_{L^{2}(I;L^{2}(Q))}^{2}\big\}\cdot\L^{d}(Q_{n,r}^{\eps})\,,
\end{align*}
where 
\[
Q_{n,r}^{\eps}:=\big\{ x\in Q\,\vert\,\dist(x,\del Q)\leq\eps nr\big\}\,.
\]
We observe that $\L^{d}(Q_{n,r}^{\eps})\to0$ as $\eps\to0$, i.e.
$\I_{Q_{n,r}^{\eps}}\to0$ point-wise $\L^{d}$-almost everywhere.
We know that $h(u^{\eps})\to h(u)$ strongly in $L^{2}(I;L^{2}(Q))$,
so dominated convergence yields that the second summand also converges
to $0$. \\
The third summand follows from two-scale convergence, i.e. $h(u_{n})\xrightharpoonup{2s,\,\mu_{\x}}h(u_{n})$
for almost every $\x$.
\end{proof}

\section{Proof of main theorem (Theorem \ref{thm:Main-Hom-Thm})\label{sec:Main-Theorem}}

\thmref{Main-Hom-Thm} is a consequence of the following.
\begin{thm}[Main theorem: homogenized limit of admissible point processes]
\label{thm:Main-Theorem}\ \\
Let $\X$ be a stationary ergodic admissible point process with distribution
$\P$ such that $\bm\X^{\complement}$ is statistically connected.
Under \assuref{PDE-Parameters}, let $u_{n}\in L^{2}(I;W^{1,2}(Q))$
be a homogenized solution from \thmref{Homogenized-System-for-n}
for the thinned point process $\X^{(n)}$.\\
Then, a subsequence of $(u_{n})_{n\in\N}$ converges to a $u\in L^{2}(I;W^{1,2}(Q))$
that is a (not necessarily unique) weak solution to the initial value
problem
\begin{align*}
\P(\bG)\del_{t}u-\nabla\cdot\left(A(u)\A\,\nabla u\right)-\lambda(\mu_{\X})h(u) & =\P(\bG)f &  & \text{in }I\times Q\,,\\
A(u)\A\nabla u\cdot\nu & =0 &  & \text{on }I\times\del Q\,,\\
u(0,x) & =\P(\bG)u_{0}(x) &  & \text{in }Q\,.
\end{align*}
Here $\A$ is the effective conductivity defined in \defref{Effective-Conductivity}
based on the event $\O=\bG=\big\{\x\in\S(\R^{d})\,\vert\,o\notin\bm\x\big\}$,
$\Omega=\S(\Rd)$, $\p=\P$ and $\lambda(\mu_{\X})$ is the intensity
of $\mu_{\X}:=\H_{\llcorner\del\bm\X}^{d-1}$.\\
Furthermore, with $\alpha_{\A}>0$ being the smallest eigenvalue of
$\A$ and $L_{h}$ being the Lipschitz constant of $h$
\begin{align*}
\esssup_{t\in I}\bign u(t)\bign_{L^{2}(Q)}^{2} & \leq\exp(C_{1})\big[\bign u_{0}\bign_{L^{2}(Q)}^{2}+C_{2}\big]\\
\bign\nabla u\bign_{L^{2}(I;L^{2}(Q))}^{2} & \leq\frac{\P(\bG)}{2\alpha_{\A}\inf(A)}\big(1+\exp(C_{1})\big)\big[\bign u_{0}\bign_{L^{2}(Q)}^{2}+C_{2}\big]\,,
\end{align*}
where
\[
C_{1}:=T\big(1+\frac{\lambda(\mu_{\X})}{\P(\bG)}(1+2L_{h})\big)\qquad\text{and}\qquad C_{2}:=\bign f\bign_{L^{2}(I;L^{2}(Q))}^{2}+2T\frac{\lambda(\mu_{\X})}{\P(\bG)}\biga h(0)\biga^{2}.
\]
\end{thm}

\begin{proof}
We note that $\A^{(n)}$ from \thmref{Homogenized-System-for-n} is
defined with cell solutions on $\Omega=\th(\Rd)$ and the push-forward
measure $\P\circ\th^{-1}$. We use the pull-back result from \thmref{pull-back-cell-solution}
to obtain a representation of $\A^{(n)}$ in terms of $\Omega=\S(\Rd)$
and the original probability distribution. 

\lemref{Convergence-of-Intensities}, \thmref{Approximation-Properties}
and \corref{Convergence-of-Effective-Conductivities} yield respectively
\[
\lambda(\mu_{\X^{(n)}})\to\lambda(\mu_{\X})\,,\qquad\P(\bG_{n})\to\P(\bG)>0\,,\qquad\A^{(n)}\to\A\,,\qquad\alpha_{\A^{(n)}}\to\alpha_{\A}>0
\]
for $\bG_{n}:=\big\{\x\,\vert\,o\notin\bm\x^{(n)}\big\}$ and $\bG:=\big\{\x\,\vert\,o\notin\bm\x\big\}$.
From the a priori estimates in \thmref{Homogenized-System-for-n},
we furthermore find 
\begin{align*}
\limsup_{n\to\infty}\left(\|u_{n}(t)\|_{L^{\infty}(0,T;L^{2}(Q))}^{2}+\|\nabla u_{n}\|_{L^{2}(I;L^{2}(Q))}^{2}\right) & <\infty\,,\\
\limsup_{n\to\infty}\bign\del_{t}u_{n}\bign_{L^{2}(I;W^{1,2}(Q)^{*})} & <\infty\,,
\end{align*}
and Aubin--Lions (or more general \thmref{Simon-Compactness}) yields
pre-compactness. These uniform bounds together with compactness arguments
yield the existence of $u\in L^{2}(I;W^{1,2}(Q))$ with generalized
time derivative $\del_{t}u\in L^{2}(I;W^{1,2}(Q))^{*}$ such that
for a subsequence
\begin{align*}
u_{n} & \xrightharpoonup[n\to\infty]{L^{2}(I;W^{1,2}(Q))}u\,, & \del_{t}u_{n} & \xrightharpoonup[n\to\infty]{L^{2}(I;W^{1,2}(Q))^{*}}\del_{t}u\,,\\
u_{n} & \xrightarrow[n\to\infty]{L^{2}(I;L^{2}(Q))}u & h(u_{n}) & \xrightarrow[n\to\infty]{L^{2}(I;L^{2}(Q))}h(u)
\end{align*}
as well as
\[
A(u_{n})\A^{(n)}\nabla u_{n}\xrightharpoonup[n\to\infty]{L^{2}(I;L^{2}(Q))}A(u)\A\nabla u\,.
\]
From here we conclude.
\end{proof}

\section{Criterion for non-degeneracy of effective conductivity \label{sec:Sufficient-Condition-for-Effective-Conductivity}}

In this chapter, we will establish a criterion for $\bm\X{}^{\complement}$
to be statistically connected (\defref{Statistically-Connected}),
that is \thmref{Percolation-Channels-Imply-Conductivity}. To be precise,
we will show that 
\[
e_{1}^{t}\A e_{1}>0
\]
as all other directions $\eta\in\R^{d}$ can be shown analogously
via rotation. The procedure will be based on \cite[Chapter 9]{kozlov1994homogenization}.
The matrix $\A$ corresponds to the matrix $\A^{0}$ there. We will
also see that $\bm\X^{\complement}$ is statistically connected iff
$\Xi\X^{\complement}$ is statistically connected.
\begin{notation*}
Given a fixed admissible point process $\X$, we write in this section
\[
\bm:=\bm\X\qquad\Xi:=\Xi\X\,.
\]
\end{notation*}
Most arguments work for more general random perforations $\Xi$ and
their filled-up versions as long as $\Xi$ has no infinite connected
component (\thmref{Filling-Holes-Preserves-Conductivity} needs additionally
that almost surely, the bounded connected components of $\R^{d}\backslash\Xi$
have non-zero distance to the infinite connected components). We refrain
from doing so since we would need to introduce the notion of stationary
random sets and the main focus here lies on point processes. 

\subsection{Variational formulation}

The following theorem gives us a different point of view on the effective
conductivity $\A$:
\begin{thm}[{Variational formulation \cite[Theorem 9.1]{kozlov1994homogenization}}]
\label{thm:Variational-Formulation-Effective-Conductivity}\ \\
For every ergodic admissible point process we have almost surely and
for every $\eta\in\R^{d}$:
\begin{align*}
\eta^{t}\A\eta=\lim_{n\to\infty}n^{-d}\inf_{v\in C_{0}^{\infty}([0,n]^{d})}\int_{[0,n]^{d}\backslash\Xi}\biga\eta-\nabla v\biga^{2}\dx x\,
\end{align*}
where $\A$ is the effective conductivity based on the event $\{o\notin\Xi\}$.
\end{thm}

\noindent The first observation we can make is that the effective
conductivity depends monotonously on the domain: The larger the set
of holes, the lower the effective conductivity. The question arises
in which cases this term becomes 0. This should only happen if $\R^{d}\backslash\Xi$
is ``insufficiently connected''. Intuitively, we want $v\approx-\eta\cdot x+\text{const}$
but at the same time, $v$ needs to be $0$ at the boundaries. If
our region is badly connected, we can hide large gradients inside
the holes, see, e.g., Figure \ref{fig:Conductivity}. As in \cite{kozlov1994homogenization},
we will see that the existence of sufficiently many ``channels''
connecting the left to the right side of a box $[0,n]^{d}$ will ensure
$e_{1}^{t}\A e_{1}>0$. Before we do that, we establish an important
fact:

\noindent 
\begin{figure}
\begin{tabular}{>{\centering}p{0.45\textwidth}>{\centering}p{0.45\textwidth}}
\includegraphics[width=0.25\paperwidth]{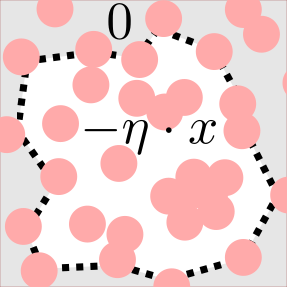} & \includegraphics[width=0.25\paperwidth]{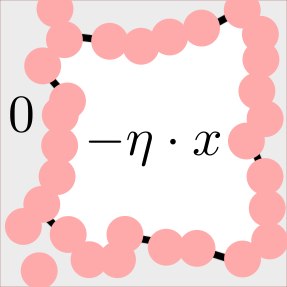}\tabularnewline
\end{tabular}

\caption{\label{fig:Conductivity}High vs low conductivity. The balls represent
$\Xi$. The white area corresponds to $v\approx-\eta\cdot x+const$.
Black lines indicate large contributions to $\int_{[0,n]^{d}\backslash\Xi}\protect\biga\eta+\nabla v\protect\biga^{2}\protect\dx x$.}
\end{figure}

\noindent We have defined statistical connectedness (\defref{Statistically-Connected})
via the filled-up Boolean model $\bm$. Unfortunately, filling up
holes is non-local (depending on the size of holes) which is troublesome
on the stochastic side. However, an analogue of \cite[Lemma 9.7]{kozlov1994homogenization}
tells us that the effective conductivity of both the Boolean model
$\Xi$ and its filled-up version $\bm$ are the same.
\begin{thm}[Filling up holes preserves the effective conductivity]
\label{thm:Filling-Holes-Preserves-Conductivity}\ \\
For every ergodic admissible point process we have almost surely
\begin{align*}
\eta^{t}\A\eta & =\lim_{n\to\infty}n^{-d}\inf_{v\in C_{0}^{\infty}([0,n]^{d})}\int_{[0,n]^{d}\backslash\bm}\biga\eta-\nabla v\biga^{2}=\lim_{n\to\infty}n^{-d}\inf_{v\in C_{0}^{\infty}([0,n]^{d})}\int_{[0,n]^{d}\backslash\Xi}\biga\eta-\nabla v\biga^{2}\,.
\end{align*}
\end{thm}

\begin{proof}
As mentioned before, this is a variation of \cite[Lemma 9.7]{kozlov1994homogenization}
fitted to our purpose. Let 
\begin{itemize}
\item $K_{n}^{s}$ be the set of islands (i.e. connected components in $\R^{d}\backslash\Xi$
of finite diameter) that intersect but do not lie inside $[0,n]^{d}$
and that are encircled by a $\Xi$-cluster of size at most $s$.
\item $L_{n}^{s}$ be the set of islands that do not completely lie inside
$[0,n]^{d}$ and that are encircled by a $\Xi$-cluster of size larger
than $s$. 
\end{itemize}
All the islands in $K_{n}^{s}$ and $L_{n}^{s}$ belong to connected
components of $\R^{d}\backslash\Xi$ different from $\bm^{\complement}$
(the unique unbounded connected component). Since $\X$ is admissible,
almost surely they all have non-zero distance to $\bm^{\complement}$.
Therefore, the following infimum decomposes, with all the infima being
over $v\in C_{0}^{\infty}([0,n]^{d})$
\begin{align*}
\inf_{v}\int_{[0,n]^{d}\backslash\Xi}\biga\eta-\nabla v\biga^{2} & =\inf_{v}\int_{([0,n]^{d}\backslash\Xi)\backslash(K_{n}^{s}\cup L_{n}^{s})}\biga\eta-\nabla v\biga^{2}+\inf_{v}\int_{K_{n}^{s}\cup L_{n}^{s}}\biga\eta-\nabla v\biga^{2}\\
 & =\inf_{v}\int_{([0,n]^{d}\backslash\bm)\backslash(K_{n}^{s}\cup L_{n}^{s})}\biga\eta-\nabla v\biga^{2}+\inf_{v}\int_{K_{n}^{s}\cup L_{n}^{s}}\biga\eta-\nabla v\biga^{2}\\
 & =\inf_{v}\int_{[0,n]^{d}\backslash\bm}\biga\eta-\nabla v\biga^{2}+C\,,
\end{align*}
with 
\[
\biga C\biga\leq\biga\eta\biga\left[\L^{d}(K_{n}^{s})+\L^{d}(L_{n}^{s})\right]
\]
and where the second equality comes from the fact that filling up
islands that lie completely inside $[0,n]^{d}$ does not change the
value of the infimum. Now, the claim follows from
\begin{itemize}
\item $\L^{d}(K_{n}^{s})\sim O(n^{d-1})$ for fixed $s$, so $\lim_{n\to\infty}n^{-d}\L^{d}(K_{n}^{s})=0$
and
\item $\lim_{n\to\infty}n^{-d}\L^{d}(L_{n}^{s})=\text{density}(L^{s})$
where $L^{s}$ denotes islands encircled by clusters of size greater
than $s$. But 
\[
\bigcap_{s\in\N}L^{s}=\emptyset,
\]
so $\lim_{s\to\infty}\text{density}(L^{s})=0$. Choosing $s$ sufficiently
large finishes the proof.
\end{itemize}
\end{proof}

\subsection{Percolation channels}
\begin{defn}[Percolation channels (see Figure \ref{fig:Percolation-channels-k-scale})]
\label{def:Percolation-Channels-1}\ \\
Fix a $k_{\scale}\in\N$. We consider the lattice $\Z_{n}^{d}\subset\Z^{d}$
and the cube with corner $z=(z_{1},\dots,z_{d})\in\Z^{d}$ 
\begin{align*}
\Z_{n}^{d}: & =\Z^{d}\cap[0,n)^{d}\qquad\text{and}\qquad\K_{z}:=\bigtimes_{i=1}^{d}[z_{i},z_{i}+1]
\end{align*}
and call two vertices $z,z'$ \emph{neighbors} if their $l^{1}$-distance
is equal to $1$.\\
We call $z$ open iff
\[
\Xi\cap k_{\scale}^{-1}\K_{z}=\emptyset.
\]
An open left-right crossing $\gamma=(z^{(1)},\dots,z^{(l)})$ of $\Z_{n}^{d}$
is called a \emph{percolation channel} in $\Z_{n}^{d}$, i.e.
\begin{enumerate}
\item all the $z^{(i)}$ are open and
\item $z_{1}^{(1)}=0$ and $z_{1}^{(l)}=n-1$.
\end{enumerate}
We define the quantity (depending on the random $\Xi$ and on $k_{\scale}$)
\begin{align*}
\mathbf{N}(n): & =\max\big\{ j\,\vert\,\gamma_{1},\dots\gamma_{j}\text{ are disjoint percolation channels in }\Z_{n}^{d}\big\}\\
 & =\text{"maximal number of disjoint percolation channels in \ensuremath{\Z_{n}^{d}}"}
\end{align*}
and the tube $L(\gamma)$ corresponding to the path $\gamma=(z^{(1)},\dots,z^{(l)})$
as
\[
L(\gamma):=\bigcup_{i\leq l}k_{\scale}^{-1}\K_{z^{(i)}}.
\]
\begin{figure}
\centering{}%
\begin{tabular}{cc}
\includegraphics[width=0.4\columnwidth]{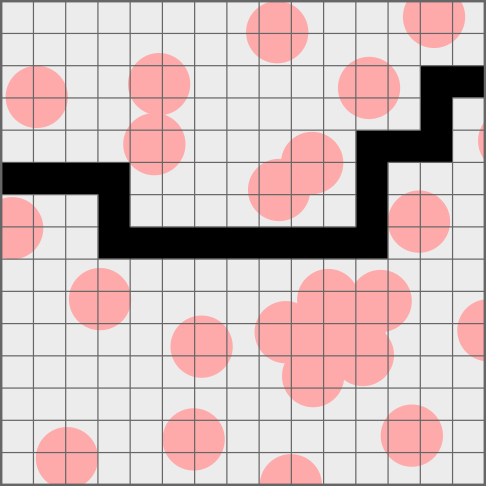} & \includegraphics[width=0.4\columnwidth]{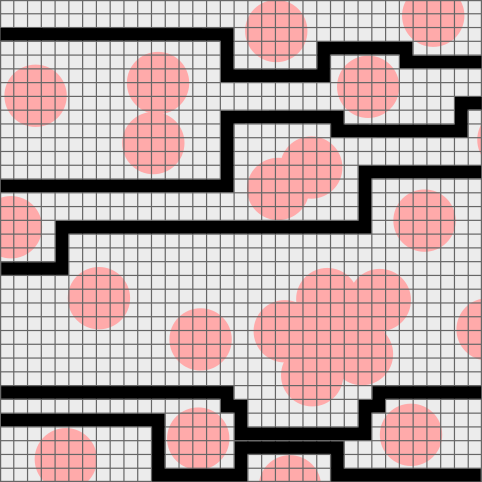}\tabularnewline
\end{tabular}\caption{\label{fig:Percolation-channels-k-scale}Percolation channels for
different $k_{\protect\scale}$}
\end{figure}
 
\end{defn}

\medskip{}
Statistical connectedness of $\R^{d}\backslash\Xi$ then reads as
follows:
\begin{thm}[Percolation channels imply conductivity]
\label{thm:Percolation-Channels-Imply-Conductivity}\ \\
For almost every realization $\x$ of an ergodic admissible point
process, we have for $\Xi=\Xi\x$
\[
\lim_{n\to\infty}n^{-d}\inf_{v\in C_{0}^{\infty}([0,n]^{d})}\int_{[0,n]^{d}\backslash\Xi}\biga e_{1}-\nabla v\biga^{2}\dx x\geq\limsup_{n\to\infty}\Big(\frac{\mathbf{N}(n)}{n^{d-1}}\Big)^{2}.
\]
In particular, the effective conductivity is strictly positive if
almost surely
\begin{equation}
\limsup_{n\to\infty}\frac{\mathbf{N}(n)}{n^{d-1}}>0\,.\label{eq:Number-of-Channels}
\end{equation}
\end{thm}

\begin{proof}
This is an analogue of \cite[Theorem 9.11]{kozlov1994homogenization}
and relies on defining a suitable vector field $\overrightarrow{F}_{\gamma}:\ [0,k_{\scale}^{-1}n]^{d}\to\R^{d}$
inside channels $\gamma=(z^{(1)},\dots,z^{(l)})$ on $\Z_{n}^{d}$.
We want $\overrightarrow{F}_{\gamma}$ to satisfy the following
\begin{itemize}
\item $\biga\overrightarrow{F}_{\gamma}(x)\biga=1$ for every $x$ inside
the tube $L(\gamma)$ and $\overrightarrow{F}_{\gamma}(x)=0$ outside.
\item $\overrightarrow{F}_{\gamma}$ is orthogonal to $\del L(\gamma)$
except on $\del L(\gamma)_{-}:=\K_{z^{(1)}}\cap\{x_{1}=0\}$ and $\del L(\gamma)_{+}:=\K_{z^{(l)}}\cap\{x_{1}=k_{\scale}^{-1}n\}$.
\item $\overrightarrow{F}_{\gamma}(x)=e_{1}$ for $x\in\del L(\gamma)_{-}\cup\del L(\gamma)_{+}$.
\item For the standard normal vector $\nu$ to $\del L(\gamma)$: 
\[
\int_{L(\gamma)}(e_{1}-\nabla v)\cdot\overrightarrow{F}_{\gamma}\dx x=\int_{\del L(\gamma)}(x_{1}-v)\overrightarrow{F}_{\gamma}\cdot\nu\dx x.
\]
\end{itemize}
Figure \ref{fig:channel} illustrates how $\overrightarrow{F}_{\gamma}$
can be chosen to satisfy these properties. 

\noindent The rest is simple. Take $\gamma_{1},\dots\gamma_{\mathbf{N}(n)}$
disjoint non-self-intersecting channels in $\Z_{n}^{d}$. Set
\[
T:=\bigcup_{i\leq\mathbf{N}(n)}L(\gamma_{i})\subset[0,k_{\scale}^{-1}n]^{d},\qquad\overrightarrow{F}:=\sum_{i\leq\mathbf{N}(n)}\overrightarrow{F}_{\gamma_{i}}.
\]
Then, 
\begin{align*}
\int_{[0,k_{\scale}^{-1}n]^{d}\backslash\Xi}\biga e_{1}-\nabla v\biga^{2}\dx x & \geq\int_{T}\biga e_{1}-\nabla v\biga^{2}\dx x\geq\int_{T}\biga(e_{1}-\nabla v)\cdot\overrightarrow{F}\biga^{2}\dx x\\
 & \geq\frac{1}{\L^{d}(T)}\left(\int_{T}(e_{1}-\nabla v)\cdot\overrightarrow{F}\dx x\right)^{2}\geq\frac{k_{\scale}^{d}}{n^{d}}\left(\int_{T}(e_{1}-\nabla v)\cdot\overrightarrow{F}\dx x\right)^{2}\,.
\end{align*}
For a fixed tube $L=L(\gamma_{i})$, we have
\begin{align*}
\int_{L}(e_{1}-\nabla v)\cdot\overrightarrow{F}\dx x & =\int_{\del L}(x_{1}-v)\overrightarrow{F}\cdot\nu\dx\H^{d-1}(x)=\int_{\del L_{-}\cup\del L_{+}}(x_{1}-v)\overrightarrow{F}\cdot\nu\dx\H^{d-1}(x)\\
 & =\int_{\del L_{+}}k_{\scale}^{-1}ne_{1}\cdot e_{1}\dx\H^{d-1}(x)=k_{\scale}^{-1}n\H^{d-1}(\del L_{+})=k_{\scale}^{-d}n.
\end{align*}
\begin{figure}
\begin{centering}
\includegraphics[width=0.8\columnwidth]{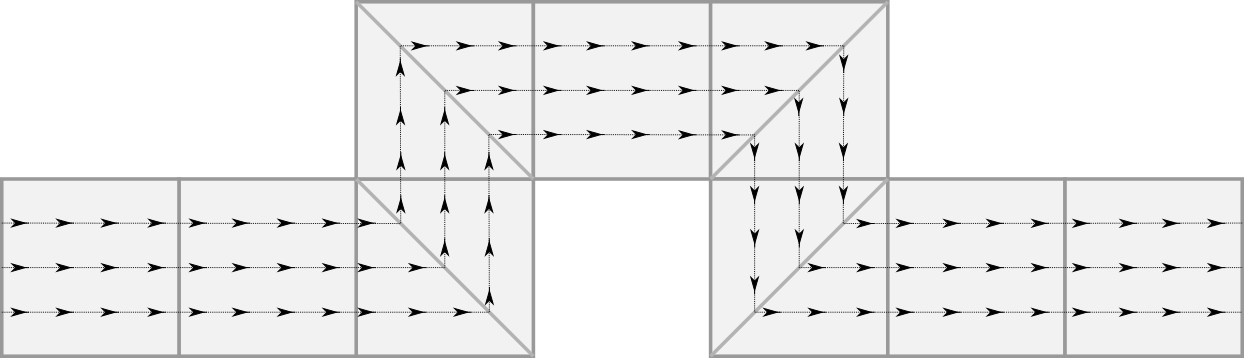}
\par\end{centering}
\caption{\label{fig:channel} Using the Gauß divergence theorem on the parts
where the vector field is constant only yields contributions from
the \textquotedblleft starting\textquotedblright{} surface and the
\textquotedblleft ending\textquotedblright{} surface. Whenever the
tube makes a left or right turn, we see that the contributions on
the diagonal surface cancel out (we have a positive contribution from
the incoming part of the tube and a negative contribution from the
outgoing one).}
\end{figure}
Therefore,
\begin{align*}
\int_{[0,k_{\scale}^{-1}n]^{d}\backslash\Xi)}\biga e_{1}-\nabla v\biga^{2}\dx x & \geq\frac{k_{\scale}^{d}}{n^{d}}\big(k_{\scale}^{-d}n\mathbf{N}(n)\big)^{2}
\end{align*}
and so
\[
\big(k_{\scale}^{-1}n\big)^{-d}\int_{[0,k_{\scale}^{-1}n]^{d}\backslash\Xi}\biga e_{1}-\nabla v\biga^{2}\dx x\geq\left(\frac{\mathbf{N}(n)}{n^{d-1}}\right)^{2}.
\]
Passing to the $\limsup$ finishes the proof.
\end{proof}
\begin{rem}[$d=2$ and bottom-top crossings]
\ \\
Let $\mathbf{L}(n)$ be the minimal number of open vertices that a
$l^{\infty}$-bottom-top crossing of $\Z_{n}^{2}$ must have. It turns
out that in $d=2$
\[
\mathbf{L}(n)=\mathbf{N}(n)\,
\]
(see \lemref{Percolation-Channels-Equals-Vertical-Crossings}). We
will use this to show \Eqref{Number-of-Channels} for the Poisson
point process $\X_{\poi}$.
\end{rem}

\section{Example: Poisson point processes \label{sec:PPP}}

The driving force behind this work has been a stationary Poisson point
process $\X_{\poi}$. It is known that the Poisson point process is
ergodic (even mixing) and its high spatial independence makes it \emph{the}
canonical random point process. As pointed out before though, $\Xi\X_{\poi}$
gives rise to numerous analytical issues which prevent the usage of
the usual homogenization tools. \\
The main theorem (\thmref{Main-Theorem}) tells us that homogenization
is still reasonable for highly irregular filled-up Boolean models
$\bm\X$ driven by admissible point processes $\X$.

\medskip{}
It is known for $\X_{\poi}$ that there exists some critical radius
$r_{c}:=r_{c}[\lambda(\X_{\poi})]\in(0,\,\infty)$ such that 
\begin{itemize}
\item $\bm\X_{\poi}$ only consists of finite clusters for $r<r_{c}$ (subcritical
regime) and
\item $\bm\X_{\poi}$ has a unique infinite cluster for $r>r_{c}$ (supercritical
regime). 
\end{itemize}
The behavior at criticality $r=r_{c}$ is still a point of research.
For details, we refer to \cite{last2017lectures} for the Poisson
point process $\X_{\poi}$ and \cite[Chapter 3]{meester1996continuum}
for the Boolean model $\Xi\X_{\poi}$.

\medskip{}
We will see in the subcritical regime that
\begin{enumerate}
\item $\X_{\poi}$ is an ergodic admissible point process and
\item $\Xi\X_{\poi}{}^{\complement}$ is statistically connected which is
equivalent to $\bm\X_{\poi}{}^{\complement}$ being statistically
connected (see \thmref{Filling-Holes-Preserves-Conductivity}).
\end{enumerate}
We therefore make the following assumption for the rest of this section:
\begin{assumption}[subcritical regime]
\label{assu:Subcritical-Regime}\ \\
We assume that 
\[
r<r_{c}\,.
\]
\end{assumption}

\begin{rem}[Scaling relation]
\ \\
$r_{c}$ has the following scaling relation 
\[
r_{c}[k^{d}\cdot\lambda(\X_{\poi})]=r_{c}[\lambda(k^{-1}\X_{\poi})]=k^{-1}r_{c}[\lambda(\X_{\poi})]\,.
\]
\end{rem}

\subsection{Admissibility of Poisson point processes}

The Mecke--Slivnyak theorem tells us that the Palm probability measure
(\thmref{-Palm-Theorem}) of a stationary Poisson point process is
just a Poisson point process with a point added in the origin. This
gives us the following lemma:
\begin{lem}[Equidistance property]
\label{cor:Equidistance-PPP}\ \\
The stationary Poisson point process $\X_{\poi}$ satisfies the equidistance
property for arbitrary $r>0$, i.e. 
\[
\P\of{\exists x,y\in\X_{\poi}\,\vert\,d(x,y)=2r}=0\,.
\]
\end{lem}

\begin{proof}
This follows from using the Palm theorem (\thmref{-Palm-Theorem})
on
\[
f(x,\x):=\sum_{x_{i}\in\x}\I\left\{ d(x,x_{i})=2r\right\} 
\]
and the Mecke--Slivnyak theorem (\cite[Theorem 9.4]{last2017lectures}).
\end{proof}

\begin{cor}[$\X_{\poi}$ is admissible]
\ \\
Under \assuref{Subcritical-Regime}, $\X_{\poi}$ is an admissible
point process.
\end{cor}

\begin{proof}
$\X_{\poi}$ is not just ergodic, but even mixing (\cite[Theorem 8.13]{last2017lectures}).
The equidistance property has been proven in \corref{Equidistance-PPP}.
Finiteness of clusters follows from the subcritical regime (\assuref{Subcritical-Regime}).
\end{proof}

\subsection{Statistical connectedness for Poisson point processes}

Proving the statistical connectedness of $\Xi\X_{\poi}{}^{\complement}$
(\defref{Statistically-Connected}) is much harder and does not immediately
follow from readily available results. Our procedure is as follows:
\begin{enumerate}
\item We employ the criterion from Section \ref{sec:Sufficient-Condition-for-Effective-Conductivity}.
Therefore, we will check that there are sufficiently many percolation
channels for $\Xi\X_{\poi}$.
\item Using the spatial independence of the Poisson point process $\X_{\poi}$,
we show that it is sufficient to only consider $2$-dimensional slices.
\item We show the statement in $d=2$ using ideas in \cite[Chapter 11]{kesten1982percolation}.
There, the result has been proven for certain iid fields on planar
graphs, including $\Z^{2}$.
\end{enumerate}
Additionally to \assuref{Subcritical-Regime}, we need sufficient
discretization for the percolation channels:
\begin{assumption}[Sufficient scaling]
\label{assu:Sufficient-Scaling}\ \\
Let $k_{\scale}\in\N$ large enough such that for the critical radius
$r_{c}$
\[
\frac{1}{2}(r_{c}-r)>\sqrt{d}k_{\scale}^{-1}\,,
\]
e.g. $k_{\scale}:=\lceil2\frac{\sqrt{d}}{r_{c}-r}\rceil+1$.
\end{assumption}

\begin{defn}[Recap and random field $(X_{z})_{z\in\Z^{d}}$]
\label{def:Random-Field-X}\ \\
Recall \defref{Percolation-Channels-1}, most importantly
\begin{align*}
\Z_{n}^{i}: & =\Z^{i}\cap[0,n)^{i}\qquad\text{and}\qquad\K_{z}:=\bigtimes_{i=1}^{d}[z_{i},z_{i}+1]
\end{align*}
as well as the notion of percolation channels for $k_{\scale}$ and
\begin{align*}
\mathbf{N}(n): & =\text{"maximal number of disjoint percolation channels in \ensuremath{\Z_{n}^{d}}"\,.}
\end{align*}
We define the \emph{random field}
\[
(X_{z})_{z\in\Z^{d}}:=\Big(\I\{\Xi\X_{\poi}\cap k_{\scale}^{-1}\K_{z}=\emptyset\}\Big){}_{z\in\Z^{d}}.
\]
We say that $z\in\Z^{d}$ is \emph{blocked} iff $X_{z}=0$ and \emph{open}
iff $X_{z}=1$. (This is consistent with \defref{Percolation-Channels-1}.)

\medskip{}
\end{defn}

\begin{thm}[Percolation channels of the Poisson point process]
\label{thm:PPP-Percolation-Channels}\ \\
Under \assuref{Subcritical-Regime} and \assuref{Sufficient-Scaling},
there is a $C>0$ such that \Eqref{Number-of-Channels} holds, i.e.
\begin{align*}
\P\big(\limsup_{n\to\infty}\frac{\mathbf{N}(n)}{n^{d-1}}\geq C\big)=1\,.
\end{align*}
In particular, $\bm\X_{\poi}{}^{\complement}$ is statistically connected
(see \thmref{Percolation-Channels-Imply-Conductivity}).
\end{thm}

The rest of the section deals with the proof of \thmref{PPP-Percolation-Channels}.
It will follow as a direct consequence of \thmref{2D-Percolation-Channels-Give-3D}
(reduction to $d=2$) and \thmref{Vertical-Crossings-Probability-Estimate}
(main result for $d=2$) which are given later.

\subsubsection{Spatial independence and moving to $d=2$}

For disjoint $U_{1},\,U_{2},\dots\,\subset\R^{d}$ and events $A_{i}$
only depending on $\X_{\poi}$ inside $U_{i}$, we know that $(A_{i})_{i}$
is an independent family. This is one of the striking properties of
a Poisson point process and we will heavily make use of it. The Boolean
model $\Xi\X_{\poi}$ for radius $r$ still retains this property
in a slightly weaker form and correspondingly the random field $(X_{z})_{z\in\Z^{d}}$: 
\begin{lem}[Independence in large distances]
\label{lem:Independence-in-Large-Distances}\ \\
Let $A,B\subset\Z^{d}$ such that
\begin{equation}
d^{\infty}(A,B):=\min_{z_{a}\in A,z_{b}\in B}\bign z_{b}-z_{a}\bign_{\infty}\geq2rk_{\scale}+1\,.\label{eq:Distance-For-Independence}
\end{equation}
Then, $(X_{z})_{z\in A}$ and $(X_{z})_{z\in B}$ are independent.
\end{lem}

\begin{proof}
$(X_{z})_{z\in A}$ is only affected by points of $\X_{\poi}$ inside
\[
U_{A}:=\bigcup_{z\in A}\Ball r{k_{\scale}^{-1}\K_{z}}\,.
\]
The same holds for $(X_{z})_{z\in B}$ and we check that \Eqref{Distance-For-Independence}
implies $U_{A}\cap U_{B}=\emptyset$.
\end{proof}
\begin{thm}[$2$-dimensional percolation channels imply channel property for $d>2$]
\label{thm:2D-Percolation-Channels-Give-3D}\ \\
 For $\tilde{z}\in\Z^{d-2}$, we define (compare to \defref{Percolation-Channels-1})
\begin{align*}
\mathbf{N}_{\tilde{z}}^{(2)}(n): & =\text{"maximal number of disjoint percolation channels in \ensuremath{\Z_{n}^{2}\times\tilde{z}}"\,.}
\end{align*}
If there are $\tilde{C},\,p_{0}>0$ such that for some $\tilde{z}\in\Z^{d-2}$
\begin{equation}
\limsup_{n\to\infty}\P\of{\mathbf{N}_{\tilde{z}}^{(2)}(n)\geq\tilde{C}n}>p_{0}>0\,,\label{eq:Number-Percolation-Channels-D=00003D2}
\end{equation}
then there exists a $C>0$ such that
\[
\limsup_{n\to\infty}\P\of{\mathbf{N}(n)\geq Cn^{d-1}}=\P\big(\limsup_{n\to\infty}\frac{\mathbf{N}(n)}{n^{d-1}}\geq C>0\big)=1\,.
\]
(This proof heavily relies on the independence structure of the Poisson
point process, i.e. \lemref{Independence-in-Large-Distances}.)
\end{thm}

\begin{proof}
$\X_{\poi}$ is stationary, so for distinct $\tilde{z}_{1},\,\tilde{z}_{2}\in\Z^{d-2}$
\[
p(n):=\P\of{\mathbf{N}_{\tilde{z}_{1}}^{(2)}(n)\geq\tilde{C}n}=\P\of{\mathbf{N}_{\tilde{z}_{2}}^{(2)}(n)\geq\tilde{C}n}.
\]
Let $k:=\lceil2rk_{\scale}\rceil+1$. By \lemref{Independence-in-Large-Distances},
the events on $\Z^{2}\times(k\tilde{z}_{1})$ are independent from
the events on $\Z^{2}\times(k\tilde{z}_{2})$. Therefore, $\big(\I\{\mathbf{N}_{k\tilde{z}}^{(2)}(kn)\geq\tilde{C}n\}\big)_{\tilde{z}\in\Z^{d-2}}$
is an iid family of Bernoulli random variables with parameter $p(n)$.
Then,
\begin{align*}
\P\big(\mathbf{N}(kn) & \geq\frac{\tilde{C}p_{0}}{2k^{d-2}}(kn)^{d-1}\big)\\
\geq & \P\big(\text{For at least \ensuremath{\frac{1}{2}p_{0}} of the }\tilde{z}\in\Z_{n}^{(d-2)}:\ \mathbf{N}_{\tilde{z}}^{(2)}(kn)\geq\tilde{C}kn\big)\\
= & \P\Big(\frac{1}{\#\Z_{n}^{(d-2)}}\sum_{\tilde{z}\in\Z_{n}^{(d-2)}}\I\{\mathbf{N}_{k\tilde{z}}^{(2)}(kn)\geq\tilde{C}n\}\geq\frac{1}{2}p_{0}\Big)\,.
\end{align*}
By \Eqref{Number-Percolation-Channels-D=00003D2} and the law of large
numbers, we get
\[
\limsup_{n\to\infty}\P\big(\mathbf{N}(n)\geq\frac{\tilde{C}p_{0}}{2k^{d-2}}n{}^{d-1}\big)\geq\limsup_{n\to\infty}\P\big(\mathbf{N}(kn)\geq\frac{\tilde{C}p_{0}}{2k^{d-2}}(kn)^{d-1}\big)=1\,.
\]
Setting $C=\frac{\tilde{C}p_{0}}{2k^{d-2}}$, we obtain \Eqref{Number-of-Channels}
after checking 
\[
\P\big(\limsup_{n\to\infty}\frac{\mathbf{N}(n)}{n^{d-1}}\geq C\big)=\limsup_{n\to\infty}\P\big(\frac{\mathbf{N}(n)}{n^{d-1}}\geq C\big)=1
\]
which finishes the proof.
\end{proof}
\begin{rem}
Spatial independence is needed to move from $d=2$ to $d\geq3$. The
strong independence properties of $\X_{\poi}$ allow far weaker conditions
on $\mathbf{N}^{(2)}(n)$ (positive probability) than on $\mathbf{N}(n)$
(probability 1). Either way, \thmref{Vertical-Crossings-Probability-Estimate}
shows that $\P(\mathbf{N}^{(2)}(n)<Cn)$ drops exponentially in $n$.
\end{rem}

\subsubsection{$d=2$: Definitions and preliminary results}

As shown before, we may limit ourselves to a fixed lattice $\Z^{2}\times0_{\Z^{d-2}}\simeq\Z^{2}$.
Therefore, we will often suppress the ``anchor point'' $0_{\Z^{d-2}}$
and just act like we are in $\Z^{2}$. Our random field from \defref{Random-Field-X}
is then by abuse of notation
\[
(X_{z})_{z\in\Z^{2}}\simeq(X_{z})_{z\in\Z^{2}\times0_{\Z^{d-2}}}\,.
\]

\begin{defn}[Vertical crossings]
\ \\
Consider the $(\Z^{2},l^{\infty})$-lattice, that is $z,z'$ are \emph{neighbors}
iff $\bign z-z'\bign_{\infty}=1$. \\
An $l^{\infty}$-bottom-top crossing in $\Z_{n}^{2}$ is called a
\emph{vertical crossing.} We call a path \emph{blocked} iff all its
vertices are blocked. We define the quantity
\[
\mathbf{L}(n):=\text{"minimal number of open vertices in a vertical crossing in \ensuremath{\Z_{n}^{2}}"}.
\]
(The percolation channels lie on the $l^{1}$-graph, while the vertical
crossings lie on the $l^{\infty}$-graph.)
\end{defn}

\begin{figure}
\begin{centering}
\begin{tabular}{cc}
\includegraphics[width=0.3\paperwidth]{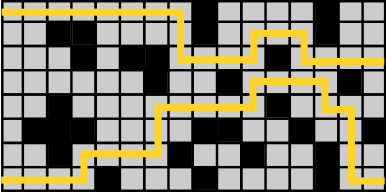} & \includegraphics[width=0.3\paperwidth]{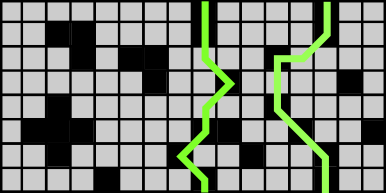}\tabularnewline
\end{tabular}
\par\end{centering}
\caption{\label{fig:percolation-channels-and-vertical-crossings}Disjoint percolation
channels vs.\textasciitilde vertical crossings. On the left side,
we see that we can only have at most two $l^{1}$-channels. The right
figure shows that any $l^{\infty}$-vertical crossing must contain
at least two open vertices. }
\end{figure}

\medskip{}
We may work with single vertical crossings instead of collections
of percolation channels: 
\begin{lem}[Percolation channels vs vertical crossings (see Figure \ref{fig:percolation-channels-and-vertical-crossings})]
\label{lem:Percolation-Channels-Equals-Vertical-Crossings}\ \\
It holds that
\[
\mathbf{N}(n)=\mathbf{L}(n).
\]
\end{lem}

\begin{proof}
See the proof of \cite[Theorem 11.1]{kesten1982percolation} based
on Menger's Theorem and \cite[Proposition 2.2]{kesten1982percolation}.
\end{proof}

\noindent The main work is proving the following equivalent of \cite[Proposition 11.1]{kesten1982percolation}:
\begin{thm}[Open vertices in vertical crossings]
\label{thm:Vertical-Crossings-Probability-Estimate}\ \\
Under \assuref{Subcritical-Regime} and \assuref{Sufficient-Scaling},
there are $C_{i}>0$ such that
\begin{align*}
 & \P\big(\exists o\leadsto\Z\times\{n\}\text{ with at most }C_{1}n\text{ open vertices}\big)\leq C_{2}\exp\big(-C_{3}n\big)\,,
\end{align*}
in particular
\[
\P\big(\mathbf{N}(n)\geq C_{1}n\big)\geq1-C_{2}n\exp(-C_{3}n)\,.
\]
\end{thm}

\noindent The proof relies on a reduction scheme of the path $\gamma:\ o\leadsto\Z\times\{n\}$.
We divide $\gamma$ into several segments which must either contain
an open vertex or contain a blocked path of large diameter. Since
we are in the subcritical regime, the probability of such paths decreases
exponentially in their diameter: 
\begin{lem}[Diameter of blocked paths]
\label{lem:Diameter-of-Blocked-Paths}\ \\
Let $z\in\Z^{2}$. Under \assuref{Subcritical-Regime} and $k_{\scale}\in\N$
as in \assuref{Sufficient-Scaling}, there are $C_{i}>0$ such that
\[
\P\big(\exists\text{blocked path }\gamma,\,z\in\gamma,\,\mathrm{diam}(\gamma)\geq n\big)\leq C_{1}\exp\big(-C_{2}n\big)\,,
\]
where 
\[
\mathrm{diam}(\gamma):=\max_{z_{1},z_{2}\in\gamma}\bign z_{1}-z_{2}\bign_{2}\,.
\]
\end{lem}

\begin{proof}
Consider the Boolean model for radius $R:=\frac{1}{2}(r+r_{c})<r_{c}$,
i.e. 
\[
\Xi^{(R)}\X_{\poi}:=\Ball R{\X_{\poi}}\,.
\]
Let $\gamma=(z^{(1)},\dots,z^{(l)})$ be a blocked path in $\Z^{2}$
containing $z$ with diameter $\geq n$. Since $\gamma$ is blocked
and $R-r>\sqrt{d}k_{\scale}^{-1}$ (\assuref{Sufficient-Scaling}),
we find for every $1\leq i\leq l$ some $x_{i}\in\X_{\poi}$ such
that
\begin{align*}
\K_{z^{(i)}}\cap & k_{\scale}\Ball r{x_{i}}\neq\emptyset
\end{align*}
and therefore
\[
\K_{z^{(i)}}\subset k_{\scale}\Ball R{x_{i}}\,.
\]
Connecting all the $z^{(i)}$ by a straight line, we obtain a continuous
path inside $k_{\scale}\Xi^{(R)}\X_{\poi}$. In particular, they all
belong to the same $k_{\scale}\Xi^{(R)}\X_{\poi}$-cluster. Then,
\begin{align*}
\P\big(\exists & \text{closed path }\gamma,\,z\in\gamma\text{ and }\mathrm{diam}(\gamma)\geq n\big)\\
 & \leq\P\big(z\text{ lies in a cluster in }k_{\scale}\Xi^{(R)}\X_{\poi}\text{ of diameter }\geq n\big)\\
 & \leq C_{1}\exp(-C_{2}n)
\end{align*}
since the occurrence of large clusters drops exponentially in their
diameter (\cite[Lemma 2.4]{meester1996continuum}).
\end{proof}

\subsubsection{Proof of \thmref{Vertical-Crossings-Probability-Estimate} (open
vertices in vertical crossings)}

\begin{figure}
\centering{}\includegraphics[width=0.85\columnwidth]{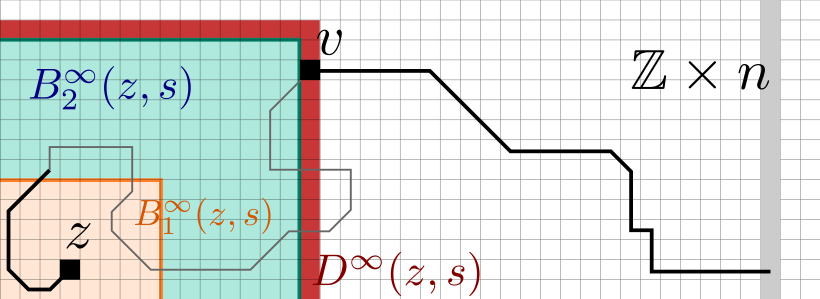}\caption{$B_{1}^{\infty},\,B_{2}^{\infty},\,D^{\infty}$ and decomposition
of paths (left-right crossing instead of top-bottom)\label{fig:Path-Decomposition}}
\end{figure}

Let $n\in\N$. As pointed out before, follow the procedure in \cite[Proposition 11.1]{kesten1982percolation}
but fitted to the continuum setting. We define $A(z,k)$ for $z\in\Z^{2}$
and $k\in\N$ as
\begin{align*}
A(z,k):=\big\{ & \exists l^{\infty}\text{-path }z\leadsto\Z\times\{n\}\text{ with at most }k\text{ open vertices}\big\}.\text{ }
\end{align*}
The idea is to break up the path $o\leadsto\Z\times\{n\}$ into multiple
segments (see \Figref{Path-Decomposition}). In each segment, we can
either reduce $k$ by $1$ or employ \lemref{Diameter-of-Blocked-Paths}.
We set
\begin{align*}
\tilde{s}: & =\lceil2rk_{\scale}\rceil+1\\
B_{1}^{\infty}(z,s): & =\big\{ v\in\Z^{2}\,\vert\,\bign z-v\bign_{\infty}\leq s\big\}\\
B_{2}^{\infty}(z,s): & =\big\{ v\in\Z^{2}\,\vert\,\bign z-v\bign_{\infty}\leq s+\tilde{s}\big\}\\
D^{\infty}(z,s): & =\big\{ v\in\Z^{2}\,\vert\,\bign z-v\bign_{\infty}=s+\tilde{s}+1\big\}="\text{boundary of }B_{2}^{\infty}(z,s)"\,.
\end{align*}
These boxes are defined so that the following holds: For fixed $z\in\Z^{2}$,
we have by \lemref{Independence-in-Large-Distances} that the random
variables $(X_{v})_{v\in B_{1}^{\infty}(z,s)}$ and $(X_{v})_{v\in\Z^{2}\backslash B_{2}^{\infty}(z,s)}$
are independent. That means the state of the vertices in $B_{1}^{\infty}(z,s)$
is independent from the state of the vertices in $\Z^{2}\backslash B_{2}^{\infty}(z,s)=B_{2}^{\infty}(z,s)^{\complement}$.
Additionally, we define the probability
\begin{align*}
g(z,s):= & \P\big(\exists z\leadsto B_{1}^{\infty}(z,s)^{\complement}\text{ blocked inside }B_{1}^{\infty}(z,s)\big).
\end{align*}
The key inequality for the iteration in $k$ is the following
\begin{equation}
\P\big(A(z,k)\big)\leq\sum_{v\in D^{\infty}(z,s)}\Big[g(z,s)\P\big(A(v,k)\big)+\P\big(A(v,k-1)\big)\Big]\label{eq:Kesten-Iteration}
\end{equation}
for $z=(z_{1},z_{2})\in\Z^{2}$ whenever $z_{2}<n-(s+\tilde{s})$.
\begin{proof}[Proof of \Eqref{Kesten-Iteration}]
 Consider the event that for some $v\in D^{\infty}(z,s)$, we find
a path $v\leadsto\Z\times\{n\}$ that has at most $k-1$ open vertices,
i.e.
\[
E:=\bigcup_{v\in D^{\infty}(z,s)}A(v,k-1)\,.
\]
Now assume that the event $A(z,k)\backslash E$ happens. Take a path
$\gamma=(z,v^{(1)},\dots,v^{(j)})$ with $v_{2}^{(j)}=n$ and at most
$k$ of the $v^{(i)}$ being open. Let $i_{1}$ be the last index
with $v^{(i_{1})}\in D^{\infty}(z,s)$. This $i_{1}$ exists since
$z_{2}<n-(s+\tilde{s})$, so $\gamma$ has to pass by $D^{\infty}(z,s)$
to reach $\Z\times\{n\}$. For this $i_{1}$, we know that $(v^{(i_{1})},\dots v^{(j)})$
completely lies in $B_{2}^{\infty}(z,s)^{\complement}$. Since $E$
does not happen, it must have $k$ open vertices. $(z,v^{(1)},\dots,v^{(i_{1})})$
is a path from $z$ to $B_{1}^{\infty}(z,s)^{\complement}$ that is
blocked everywhere except its end. Therefore,
\begin{align*}
A(z,k)\backslash E\subset\bigcup_{v\in D^{\infty}(z,s)}\big\{ & \exists z\leadsto B_{1}^{\infty}(z,s)^{\complement}\text{ blocked in }B_{1}^{\infty}(z,s)\text{ and}\\
 & \exists v\leadsto\Z\times\{n\}\text{ in }B_{2}^{\infty}(z,s)^{\complement}\text{ with }\text{at most }k\text{ open vertices}\big\}.
\end{align*}
As mentioned before, the events in $B_{1}^{\infty}(z,s)$ and $B_{2}^{\infty}(z,s)^{\complement}$
are independent from each other. This gives us
\begin{align*}
\P\big(A(z,k)\backslash E\big)\leq\sum_{v\in D^{\infty}(z,s)}\P\big( & \exists z\leadsto B_{1}^{\infty}(z,s)^{\complement}\text{ blocked in }B_{1}^{\infty}(z,s)\text{ and}\\
 & \exists v\leadsto\Z\times\{n\}\text{ in }B_{2}^{\infty}(z,s)^{\complement}\text{ with }\text{at most }k\text{ open vertices}\big)\\
=\sum_{v\in D^{\infty}(z,s)}\P\big( & \exists z\leadsto B_{1}^{\infty}(z,s)^{\complement}\text{ blocked in }B_{1}^{\infty}(z,s)\big)\\
\times\P\big( & \exists v\leadsto\Z\times\{n\}\text{ in }B_{2}^{\infty}(z,s)^{\complement}\text{ with }\text{at most }k\text{ open vertices}\big)\\
\leq\sum_{v\in D^{\infty}(z,s)}g( & z,s)\P\big(A(v,k)\big).
\end{align*}
\begin{align*}
\P\big(A(z,k)\big) & \leq\P\big(A(z,k)\backslash E\big)+\P(E)\leq\sum_{v\in D^{\infty}(z,s)}\Big[g(z,s)\P\big(A(v,k)\big)+\P\big(A(v,k-1)\big)\Big]
\end{align*}
concludes the proof of \Eqref{Kesten-Iteration}.
\end{proof}
\medskip{}

Observe that the reduction in $k$ can only happen until $k=0$, so
more $g(z,s)$-terms have to show up at some point. Since any path
$z\leadsto B_{1}^{\infty}(z,s)^{\complement}$ has diameter of at
least $s$, \lemref{Diameter-of-Blocked-Paths} tells us that $g(z,s)\leq C_{1}\exp\big(-C_{2}s\big)$
for $C_{1},\,C_{2}>0$ independent of $z$ and $s$. Choose $s$ large
such that
\begin{align*}
g(z,s)\cdot\#D^{\infty}(z,s) & \leq C_{1}\exp\big(-C_{2}s\big)\,\cdot8(s+\tilde{s}+1)\leq\frac{1}{4}.
\end{align*}
For simplicity, we introduce 
\begin{align*}
D^{\infty} & :=D^{\infty}(o,s)\qquad\text{and}\qquad h(z,y):=\begin{cases}
g(z,s) & \text{if }y=0\\
1 & \text{if }y=1.
\end{cases}
\end{align*}
and rewrite \Eqref{Kesten-Iteration} into
\begin{align}
\P\big(A(z,k)\big) & \leq\sum_{v^{(1)}\in D^{\infty},\,y_{1}\in\{0,1\}}h(z,y_{1})\P\big(A(z+v^{(1)},k-y_{1})\big)\,.\label{eq:Kesten-Iteration-2}
\end{align}
We now iteratively use \Eqref{Kesten-Iteration-2} up to $l$ times
as it is only applicable when $z_{2}<n-(s+\tilde{s})$. All the $v^{(i)}$
are summed over $D^{\infty}$ and all the $y_{i}$ over $\{0,1\}$.
\begin{align*}
\P\big(A & (o,k)\big)\\
\leq\ \  & \sum_{\substack{i\leq l,v^{(i)},y_{i}\\
v_{2}^{(1)}+\dots+v_{2}^{(l)}<n-(s+\tilde{s})\\
y_{1}+\dots+y_{l}\leq k
}
}\P(A(v^{(1)}+\dots+v^{(l)},k-y_{1}-\dots-y_{l}))\prod_{m\leq l}h(v^{(1)}+\dots+v^{(m-1)},y_{m})\\
+ & \sum_{\frac{n-(s+\tilde{s})}{s+\tilde{s}+1}\leq j\leq l}\sum_{\substack{i\leq j,v^{(i)},y_{i}\\
v_{2}^{(1)}+\dots+v_{2}^{(j-1)}<n-(s+\tilde{s})\\
v_{2}^{(1)}+\dots+v_{2}^{(j)}\geq n-(s+\tilde{s})\\
y_{1}+\dots+y_{j}\leq k
}
}\P(A(v^{(1)}+\dots+v^{(j)},k-y_{1}-\dots-y_{j}))\prod_{m\leq j}h(v^{(1)}+\dots+v^{(m-1)},y_{m})\\
\leq\ \  & \sum_{\substack{y_{1},\dots,y_{l}\\
y_{1}+\dots+y_{l}\leq k
}
}\Big(\#D^{\infty}\sup_{\substack{v\in\Z^{2}\\
v_{2}<n-(s+\tilde{s})
}
}h(v,y_{m})\Big)^{l}+\sum_{\frac{n-(s+\tilde{s})}{s+\tilde{s}+1}\leq j\leq l}\sum_{\substack{y_{1},\dots,y_{j}\\
y_{1}+\dots+y_{j}\leq k
}
}\Big(\#D^{\infty}\sup_{\substack{v\in\Z^{2}\\
v_{2}<n-(s+\tilde{s})
}
}h(v,y_{m})\Big)^{j}\\
\end{align*}
We iterate as long as $0+v_{2}^{(1)}+\dots v_{2}^{(j)}<n-(s+\tilde{s})$,
otherwise we stop for $0+v_{2}^{(1)}+\dots v_{2}^{(j)}$ and land
in the second summand. Only $y_{1}+\dots+y_{j}\leq k$ matters since
$A(z,m)=0$ whenever $m<0$. Also observe that $v_{2}^{(1)}+\dots+v_{2}^{(j)}\geq n-(s+\tilde{s})$
can only happen if 
\[
j\geq\frac{n-(s+\tilde{s})}{s+\tilde{s}+1}
\]
since we ``gain'' at most $s+\tilde{s}+1$ to the second component
in each $v^{(i)}$.\\
Let $\alpha>0$ be large enough such that
\begin{align*}
\phi(\alpha): & =\sum_{y\in\{0,1\}}\sup_{\substack{z\in\Z^{2}\\
z_{2}<n-(s+\tilde{s})
}
}\#D^{\infty}\cdot h(z,y)\cdot e^{-\alpha y}\leq\frac{1}{4}+\#D^{\infty}\cdot e^{-\alpha}\leq\frac{1}{2}.
\end{align*}
Then,
\begin{align*}
\P\big(A(o,k)\big)\leq\ \  & e^{\alpha k}\sum_{\substack{y_{1},\dots,y_{l}\\
y_{1}+\dots+y_{l}\leq k
}
}\Big(\#D^{\infty}\sup_{\substack{v\in\Z^{2}\\
v_{2}<n-(s+\tilde{s})
}
}h(v,y_{m})\Big)^{l}\prod_{i\leq l}e^{-\alpha y_{i}}\\
+ & e^{\alpha k}\sum_{\frac{n-(s+\tilde{s})}{s+\tilde{s}+1}\leq j\leq l}\sum_{\substack{y_{1},\dots,y_{j}\\
y_{1}+\dots+y_{j}\leq k
}
}\Big(\#D^{\infty}\sup_{\substack{v\in\Z^{2}\\
v_{2}<n-(s+\tilde{s})
}
}h(v,y_{m})\Big)^{j}\prod_{i\leq j}e^{-\alpha y_{i}}\\
\leq\ \  & e^{\alpha k}\phi(\alpha)^{l}+e^{\alpha k}\sum_{j\geq\frac{n-(s+\tilde{s})}{s+\tilde{s}+1}}^{l}\phi(\alpha)^{j}\leq e^{\alpha k}\big[2^{-l}+2^{-\frac{n-(s+\tilde{s})}{s+\tilde{s}+1}+1}\big]\,.
\end{align*}
Since $l\in\N$ was arbitrary, we get
\begin{align*}
\P\big(A(o,k)\big) & \leq e^{\alpha k}\cdot2^{-\frac{n-(s+\tilde{s})}{s+\tilde{s}+1}+1}=e^{\alpha k-\big[\frac{n-(s+\tilde{s})}{s+\tilde{s}+1}-1\big]\ln2}=C_{3}e^{\alpha k-C_{4}n}.
\end{align*}
Now we finally make use of $k$. Setting $C_{5}:=\frac{C_{4}}{2\alpha}$
and $k:=C_{5}n$, we obtain the claim
\[
\P\big(A(o,C_{5}n)\big)\leq C_{3}\exp\big(-\frac{1}{2}C_{4}n\big)\,.
\]

\medskip{}

\begin{acknowledgement*}
This work was funded by the German Leibniz Association via the Leibniz
Competition 2020. 
\end{acknowledgement*}

%\printbibliography[heading=bibintoc]

\bibliography{../stochhom_literature.bib}
\bibliographystyle{alpha}

\end{document}